\documentclass[11pt]{article}        
\usepackage{amsmath,amsthm}
\usepackage{amsfonts}
\usepackage{amssymb}
\usepackage{xcolor,graphics,graphicx,subfigure}

\newtheorem{theorem}{Theorem}[section]

\newtheorem{remark}[theorem]{Remark}
\usepackage[ruled,vlined]{algorithm2e}
\usepackage{booktabs,caption}

\usepackage{xcolor}
\newcommand{\rev}[1]{{#1}}
\newcommand{\revv}[1]{{#1}}

\begin{document}

\title{A Novel Stochastic Interacting Particle-Field Algorithm for 3D Parabolic-Parabolic Keller-Segel Chemotaxis System}

\author{Zhongjian Wang \thanks{   Division of Mathematical Sciences, School of Physical and Mathematical Sciences, Nanyang Technological University, 21 Nanyang Link, Singapore 637371. zhongjian.wang@ntu.edu.sg}       \and
        Jack Xin \thanks{Department of Mathematics, University of California at Irvine, Irvine, CA 92697, USA. jack.xin@uci.edu} \and
        Zhiwen Zhang \thanks{Department of Mathematics, The University of Hong Kong, Pokfulam Road, Hong Kong SAR, P.R. China.  Materials Innovation Institute for Life Sciences and Energy (MILES), HKU-SIRI, Shenzhen, P.R. China. zhangzw@hku.hk}
}

\maketitle

\begin{abstract}
We introduce an efficient stochastic interacting particle-field (SIPF) algorithm with no history dependence for computing aggregation patterns and near singular solutions of parabolic-parabolic Keller-Segel (KS) chemotaxis system in three-dimensional (3D) space. In our algorithm, the KS solutions are approximated as empirical measures of particles coupled with a smoother field (concentration of chemo-attractant) variable computed by a spectral method. Instead of using heat kernels that cause history dependence and high memory cost, we leverage the implicit Euler discretization to derive a one-step recursion in time for stochastic particle positions and the field variable based on the explicit Green's function of an elliptic operator of the form Laplacian minus a positive constant. In numerical experiments, we observe that the resulting SIPF algorithm is convergent and self-adaptive to the high-gradient part of solutions. Despite the lack of analytical knowledge (such as a self-similar ansatz) of a blowup, the SIPF algorithm provides a low-cost approach to studying the emergence of finite-time blowup in 3D space using only dozens of Fourier modes and by varying the amount of initial mass and tracking the evolution of the field variable. Notably, the algorithm can handle multi-modal initial data and the subsequent complex evolution involving the merging of particle clusters and the formation of a finite time singularity with ease.

{\textbf{Keywords:} fully parabolic Keller-Segel system \and interacting particle-field approximation \and singularity detection \and critical mass \and finite-time blowup.
}

\end{abstract}

	\section{Introduction}
	Chemotaxis partial differential equations (PDEs) were introduced by Keller and Segel (KS \cite{keller1970initiation}) to describe the aggregation of the slime mold amoeba Dictyostelium discoideum due to an attractive chemical substance. A related random walk model was known earlier, developed by Patlak \cite{Patlak_53}.  For an analysis of basic taxis behaviors (such as aggregation, blowup, and collapse) based on reinforced random walks, we refer to \cite{Oth_97}. In this work, we consider the parabolic-parabolic (fully parabolic) KS system of the following form:
	\begin{align}
		\rho_{t}  &= \nabla \cdot (\mu \, \nabla \rho - \chi\, \rho\,  \nabla c ), \nonumber \\ 
		\epsilon \,   c_t  &= \Delta\, c - k^2 \, c + \rho,  \label{ppKS}
	\end{align}
	where $\chi, \mu $ ($\epsilon, k$) are positive (non-negative) constants. The model is called elliptic if $\epsilon =0$ (when $c$ evolves rapidly to a local equilibrium), and parabolic if $\epsilon >0$. 
	Here, $\rho$ is the density of active particles (such as bacteria), and $c$ is the concentration of
	chemo-attractant (e.g. food). For a more detailed discussion, please refer to Section \ref{sec:KSmodel}.
	
	\medskip 
	
	The KS systems \eqref{ppKS} have been extensively studied for several decades, with various cases and dimensions explored. \rev{It is well known that the system \eqref{ppKS} in $\mathbb{R}^2$ converges to Dirac-delta type function in finite time given sufficiently large initial mass in the parabolic-elliptic case ($k=0$ and $\epsilon=0$) \cite{perthame2006transport} and the fully parabolic case \cite{herrero1997blow}. In addition to the $\delta$ type singularity,  $C(T-t+|x|^2)^{-1}$ is shown to be an alternative profile for self-similar blowup \cite{giga2011asymptotic,souplet2019blow} in 3D parabolic elliptic case. Given the smallness assumption of the initial data, \cite{calvez2008parabolic,mizoguchi2013global,lemarie2013small,takeuchi2021keller} established global well-posedness results for 3D fully parabolic KS in various function spaces.  To date, 
    the complete characterization (criteria and profile) of the finite time 3D singularity of (\ref{ppKS}) with $\epsilon > 0$ remains largely open.} \rev{See also  \cite{hillen2001global} for a related model that does not blow up and still exhibits spiky solutions by assuming a saturation concentration for the bacteria. }

		\medskip 
	Many notable numerical methods have been developed for KS systems. 
	Chertock et al. \cite{chertock2008second} developed a finite-volume method on a class of chemotaxis and haptotaxis models for accurate and efficient simulations. 
    Shen et al. \cite{shen2020unconditionally} proposed an energy dissipation and bound preserving scheme that is not restricted to specific spatial discretization methods. The bound preserving property is achieved through modifications of the system. 
    Chen et al. \cite{chen2022error} developed a fully-discrete finite element method (FEM) scheme for the 2D parabolic-elliptic KS, following the approach of Shen et al. \cite{shen2020unconditionally}. They showed that the proposed scheme will blow up in a finite time, under assumptions similar to those in the continuous blow-up scenarios. In the classic setting, Liu and Wang \cite{liu2018positivity} reformulated the equation using the Le Châtelier Principle to attain a positive-preserving scheme.  See \cite{DG_KS_2017,DG_KS_2024} among others on discontinuous Galerkin methods.
    Besides the aforementioned numerical methods 
    for 2D KS, \rev{we refer to 
    \cite{PPFEM3D_2013,EpXia_2019} on mesh-based (finite element, finite volume/difference with difference of potential domain decomposition) methods for 3D fully parabolic KS models.}

\medskip

	As an alternative to the Eulerian discretization methods mentioned above, there have been steady advancements in the Lagrangian formulations for the KS system (\ref{ppKS}) and related equations. \rev{
    The Lagrangian framework approximates $\rho$ in \eqref{ppKS} as the density of an $N$-interacting particle system as $N$ tends to infinity.}
	Stevens \cite{stevens2000derivation} derived a convergent $N$-particle system for the fully parabolic KS.
    \rev{Craig and Bertozzi \cite{craig2016blob} proved the convergence of a blob method for a related aggregation equation.  Liu et al. \cite{liu2019propagation,liu2017random} developed a random particle blob method with a mollified kernel for the parabolic-elliptic KS, and proved its convergence when the limiting (macroscopic mean field) equation admits a global weak solution. See also \cite{chen2022mean} for an extension to the fully parabolic KS. As noted in 
    \cite{mischler2013kac}, the success of this line of analysis relies on the limiting nonlinear mean field equation, rather than the underlying many-particle Markov process.}
\medskip

    \rev{On numerically implemented particle methods for KS systems, Havskovéč- Ševčovič \cite{havskovec2009stochastic,havskovec2011convergence} developed a convergent regularized particle system for the 2D parabolic-elliptic KS and in the presence of a passive flow \cite{chemomix_yao}. The authors of this paper conducted a deep learning study on chemotaxis aggregation in 3D laminar and chaotic flows, based on a kernel regularized particle method for parabolic-elliptic KS systems \cite{DPKS_22}. For 2D fully parabolic KS, \cite{fatkullin2012study} proposed a particle-field method and \cite{tomasevic:tel-01932777} investigated a fully particle approach.
    Though both methods may be generalized to higher dimensions,  the former suffers from numerical instability when capturing blowup behavior and the latter is inefficient for long-time simulation due to memory costs. See Section \ref{sec:lagrangian} for an in-depth discussion and comparison with our proposed method.}

	

    \medskip
	
	\rev{To the best of our knowledge, an efficient particle-based (mesh-free) method 
    capable of characterizing blow-up behaviors of 3D fully parabolic KS, especially the critical mass, is unknown.  In this paper, we propose a novel stochastic interacting particle-field (SIPF) algorithm
    to serve this purpose}. 
    \medskip
    
    Our method takes into account the coupled stochastic particle evolution (density $\rho$) and the accompanying field (concentration $c$) in the system and allows for a self-adaptive simulation of focusing and potentially singular behavior. In the SIPF algorithm, we represent the active particle density $\rho$ using empirical particles, while the concentration field $c$ is discretized using a spectral method instead of a finite difference method \cite{fatkullin2012study}. This is possible since the field $c$ is smoother than density $\rho$ \rev{and does not require interpolation of $\rho$ to grid points as \cite{fatkullin2012study}}. We demonstrate the effectiveness of our method through numerical experiments in 3D space, which, to the best of our knowledge, have not been systematically computed and benchmarked before.
	\medskip
    
	It is worth noting that pseudo-spectral methods were employed to compute the nearly singular solutions of the 3D Euler equations \cite{hou2007computing}. Subsequently, the finite-time blowup of 3D axisymmetric Euler equations was computed using the adaptive moving mesh method \cite{luo2014toward}. These methods represent the cutting edge in the computation of nearly singular solutions of 3D Euler equations. Nevertheless, we also point out that the implementation of pseudo-spectral methods for 3D problems demands substantial computational resources, while the adaptive moving mesh method requires sophisticated design and advanced programming skills.  
	\medskip
    
	It is important to recognize that the Lagrangian algorithms used in the computation of parabolic-elliptic KS systems, such as the one proposed in \cite{havskovec2009stochastic}, cannot be directly generalized to the fully parabolic KS. These algorithms rely on the assumption that the field $c$ at time $t$ can be accessed through particle density $\rho$ at the same instant. Hence only the local update of the particle density is required. However, a direct generalization to the fully parabolic KS \rev{(e.g. \cite{tomasevic:tel-01932777})}  will require historical particle density $\rho$ from the starting time of the algorithm. An example and related convergence analyses can be found in \cite{chen2022mean}. Nonetheless, from a computational perspective, the volume of such historical data grows over time, posing a costly burden on memory and computational resources. In contrast, our {\it SIPF algorithm computes particle and field values only once per time step, without relying on a long history.} Therefore,  our computational cost does not increase over time.
	
\medskip
    
	The main objective of this paper is to propose a novel stochastic interacting particle-field algorithm for the fully parabolic KS system. While we provide \rev{stability analysis and }numerical verification of the convergence of the SIPF algorithm, a detailed theoretical analysis will be left as a future work.
	\medskip
    
	The rest of the paper is organized as follows. In Section 2, we briefly review the blow-up behavior in the fully parabolic KS models under critical mass conditions and the Lagrangian formulations for computation. Section 3 outlines our proposed SIPF algorithms for solving the fully parabolic KS system. We simplify a theoretically equivalent method with history-dependent parabolic kernel functions, which is computationally undesirable, into efficient recursions. In Section 4, we present numerical results to demonstrate the performance of our method for both radial and multi-modal initial data. \rev{The SIFP results are in agreement with those from a fully resolved finite difference method in the radially symmetric case.} Concluding remarks are in Section 5.
	
	\section{Parabolic-Parabolic KS System} \label{sec:KSmodel}
	In this section, we provide some theoretical analyses of singular behaviors and related computational methods for KS models in both parabolic-elliptic cases and parabolic-parabolic (fully parabolic) cases.  We begin by revisiting the KS model:
	\begin{align}
		\rho_{t} &= \nabla \cdot (\mu \, \nabla \rho - \chi\, \rho\,  \nabla c ), \label{KS_rho} \\ 
		\epsilon   c_t  &= \Delta\, c - k^2 \, c + \rho,  \label{KS_c}\\
		x\in & \Omega \subseteq \mathbb{R}^d,\quad t\in[0,T].
	\end{align}
	The first equation \eqref{KS_rho} of $\rho$ models the evolution of the density of active particles, such as bacteria. The bacteria diffuse with mobility $\mu$ and drift in the direction of $\nabla c$ with a velocity $\chi \nabla c$, where $\chi$ is called chemo-sensitivity. The second equation \eqref{KS_c} of $c$ models the evolution of the concentration of chemo-attractant, such as food. The increment in $c$ is proportion to $\rho$, which indicates the aggregation or attraction between active particles. An additional important physical parameter in the model is $\epsilon$ in Eq.\eqref{KS_c}, which represents the timescale of the chemotaxis. When $\epsilon\not=0$, the system is referred to as a parabolic-parabolic KS system. For $\epsilon=0$ the system is reduced to the parabolic-elliptic case, which assumes that the chemical attractant released by the active particle instantaneously reaches equilibrium.
	
	\subsection{From Critical Collapse to  Coexistence of Blow-up and Global Smooth Solutions}\label{sec:blowup}
	Well-known KS dichotomy (critical collapse) states that $8\pi$ is the critical mass for the simplest two-dimensional parabolic-elliptic KS system in $\Omega=\mathbb{R}^2$, namely \eqref{ppKS} with $\epsilon=k=0$,
	
	The well-known KS dichotomy, also known as critical collapse, states that in the simplest 2D parabolic-elliptic KS system in the domain $\Omega=\mathbb{R}^2$, the critical mass is $8\pi$. This system, denoted by \eqref{ppKS} with $\epsilon=k=0$, i.e., 
	\begin{align}
		\rho_{t} &= \nabla \cdot ( \nabla \rho -  \rho\,  \nabla c ), \nonumber \\ 
		\Delta\, c  &=- \rho,  \label{ppKS2Dsimple}
	\end{align}
     exhibits a critical behavior where the solution either remains global and bounded or blows up in finite time depending on the initial mass of the density $\rho$. Specifically, we have  
	\begin{enumerate}
		\item If $M_0<8\pi$, the system has a global smooth solution.
		\item If $M_0>8\pi$, the system blows up in finite time in the sense of $|\cdot|_\infty$ norm.
	\end{enumerate}
	It can be seen from the classical variance identity for system \eqref{ppKS2Dsimple}, \cite{perthame2006transport}, that,
	\begin{align}
		\frac{d}{dt}\int_{x\in \mathbb{R}^2} |x|^2 \, \rho(x)\, dx=\frac{M}{2\pi}(8\pi-M).\label{simple_mass_result}
	\end{align}
	Then, the solution of \eqref{ppKS2Dsimple} exhibits a quantized concentration of mass at the origin, which is a type of blow-up known as a $\delta$-blow up. 
	
	In the case of the system \eqref{ppKS2Dsimple} on $\mathbb{R}^d$ with $d\geq 3$, the identity \eqref{simple_mass_result} is no longer applicable, and the evolution of the KS system is not as straightforward. However, the coexistence of blow-up and global smooth solutions still depends on the size of the initial data. 
 \rev{For fully parabolic KS, \cite{winkler2013finite,winkler_blow-up_2020} show that the system may blow up in finite time over a large set of radial initial data. On the other side of the dichotomy, it has been shown \cite{takeuchi2021keller} that global strong solutions exist for small initial data in the fully parabolic system \eqref{ppKS}.}
 
 In addition, the blow-up profile can differ from the $\delta$-type blow-up observed in 2D cases. For instance, it has been shown in \cite{herrero1998self} that in 3D parabolic elliptic systems, there exist radial, positive, backward self-similar solutions of the form,
	\begin{align}
		\rho(x,t)=\frac{V(x/\sqrt{T-t})}{T-t}, \quad \quad 0<t<T,
	\end{align}
	where the radially decreasing profile function $V$ satisfies $\lim_{y\to\infty}y^2V(y)=L\in \mathbb{R}^{+}$. 
	Later in a more refined result by \cite{souplet2019blow}, the blowup is said to be type I if 
	\begin{align}
		0 <\limsup_{t\to T}\, (T-t)\, \|\rho\|_{\infty}\, <\infty.
	\end{align}
	Then for radial initial data in $L^1(\mathbb{R}^3)$, if a blowup is type I, $\exists \, C>0$ such that 
	\begin{align}
		\rho(x,t)\leq C(T-t+|x|^2)^{-1}, \quad 0<|x|\leq R, \quad 0<t<T.
	\end{align}
	 
	\rev{To the best of our knowledge, the complete characterization of blow-up criteria with non-radial initial data and the profile of the fully parabolic KS system in $\mathbb{R}^3$ 
 has not been analyzed, also as discussed in \cite{winkler_blow-up_2020}.} Therefore, numerical computations are necessary to investigate the potential singular behavior, which we will discuss in Section \ref{sec:blowup_num}.

	\subsection{Lagrangian formulations}\label{sec:lagrangian}
	As a fundamental step of deriving the algorithms, we begin by introducing the Lagrangian formulation of the active particle density $\rho$ in the KS system \eqref{ppKS}. We focus on
	the elliptic system with $\epsilon=k=0$, specifically \eqref{ppKS2Dsimple}, which can be 
	generalized to any dimension $d$. By considering the equation $\Delta c = -\rho$ and utilizing the Green's function of the Laplacian operator in $\mathbb{R}^d$, we can deduce the following:
	\begin{align}
		c(x,t)=\left\{\begin{array}{l}
			-\frac{1}{2\pi}\int \ln|x-y| \, \rho(y,t), \quad d=2\\
			C_d\int \frac{1}{|x-y|^{d-2}}\, \rho(y,t)\, dy, \quad d\geq 3 
		\end{array}\right.,
	\end{align}
	where $C_d=\frac{\Gamma(d/2+1)}{d(d-2)\pi^{d/2}}$. So the convection term in \eqref{KS_rho} can be expressed as follows:
	\begin{align}
		\nabla c(x) = -\frac{\Gamma(d/2)}{2\pi^{d/2}}\int \frac{x-y}{|x-y|^d}\, \rho(y,t)\, dy.
	\end{align}
	Now we arrive at the interactive stochastic differential equation system of $P$ particles, $\{X^p_t\}_{p=1:P}$,
	\begin{align}
		dX^p_t=-\chi\frac{M}{P}\sum_{q\not=p} \frac{\Gamma(d/2)}{2\pi^{d/2}} \frac{X^p_t-X^q_t}{|X^p_t-X^q_t|^d}+\sqrt{2\mu} \, dW^p_t, \quad p=1,\cdots,P, \label{IPS_SDE0}
	\end{align}
	where $W^p_t$ denotes independent identically distributed standard Brownian motions.
	In \cite{mischler2013kac}, it has been demonstrated, under mild regularity condition, that as $P\to \infty$, the distribution of empirical particles $\{X^p_t\}_{p=1:P}$ converges to $\rho$ in the continuous PDE system \eqref{KS_rho}. To study the singularity behavior in the parabolic elliptic KS systems, several novel numerical methods have been developed and implemented, as discussed in  \cite{liu2017random,havskovec2009stochastic,DPKS_22}. 
	
	In the fully parabolic case ($\epsilon\not  =0$), the solution of chemical concentration $c$
	is obtained by solving a parabolic equation, which is no longer Markovian as in \eqref{IPS_SDE0}. In this case, at time $t>0$, the solution of $\rho$ in the interval $[0,t]$ needs to be 
	incorporated in the representation of $c$. Specifically, we have to consider the following: 
	\begin{align}\label{parabolic_c}
		c(\cdot,t)=e^{-k^2t}e^{t\Delta} c(\cdot,0) + \int_0^t e^{k^2(s-t)}e^{(t-s)\Delta}\, \rho(\cdot,s)\, ds,
	\end{align}
	where the heat semigroup operator $e^{t\Delta}$ is defined by
	\begin{align}
		(e^{t\Delta}f)(x,t):=\int \frac{e^{-|x-y|^2/(4t)}}{(4\pi t)^{d/2}}\, f(y)\, dy.
	\end{align}
	Similar to \eqref{IPS_SDE0}, the empirical particle system  converging to density $\rho$ reads:
	\begin{align}
		d X^{p}_t = \,  {\chi } \nabla_{X}\, c(X^{p}_{t},t) 
		\, dt 
		+ \sqrt{2\,\mu} \, d W^p, \;\; p=1,\cdots, P,  \label{IPS_SDE}
	\end{align}
	and $W^p$'s are independent Brownian motions in $\mathbb{R}^d$. However, due to the historic path dependence in the solution of $c$ in \eqref{parabolic_c}, direct computation of the drift $\nabla_{X}\, c(X^{p}_{t},t) $ in \eqref{IPS_SDE} \rev{as \eqref{IPS_SDE0} would result in significant memory and computational cost in each step, which increases with computational time $T$. More precisely, by directly substituting \eqref{parabolic_c} into \eqref{IPS_SDE}, the empirical particle system is as follows,
 \begin{align}\label{IPS_SDE_full}
     d X^{p}_t = &\,  {\chi } \,e^{-k^2t}e^{t\Delta}\nabla_{X} c(X^{p}_t,0) dt + \sqrt{2\,\mu} \, d W^p, \nonumber\\
     &-\left( \chi\frac{M}{P}\int_0^t e^{k^2(s-t)}\sum_{q\not=p} \frac{1}{(4\pi (t-s))^{d/2}}\frac{X^p_t-X^q_s}{2(t-s)}\, ds \right)		 dt  .
 \end{align}
In the discrete scheme of \eqref{IPS_SDE_full}, at the $n$-th temporal step, one shall compute the interaction of $X_{t_n}$ with other particle positions over time interval $[0,t_n]$. At the $n$-th step, the computational cost is $\mathcal{O}(n P^2)$ and memory cost is $\mathcal{O}(nP)$. } 
 
 \rev{A purely probabilistic particle method based on \eqref{IPS_SDE_full} for 2D fully parabolic KS system is proposed in \cite{tomasevic:tel-01932777}. The method is efficient due to the interaction, as mentioned by \cite{tomasevic:tel-01932777}. In contrast, \cite{fatkullin2012study} proposed a memory-less approach for 2D systems similar to ours. However, the method in \cite{fatkullin2012study} computes $c$ only on the spatial grid points, which leads to numerical inaccuracy in moving the particles. More details for computing 2D systems can be found in \cite{tomasevic:tel-01932777}.  To the best of our knowledge, a memory-less algorithm (memory cost $\mathcal{O}(P)$ for saving particle positions) to compute fully parabolic KS systems in 3D, or higher dimensions, has not been developed yet.} We will present such an algorithm in the following section.

	\section{SIPF Algorithms for  Parabolic-Parabolic KS }
	
	In this section, we will present the SIPF algorithm for solving the fully parabolic KS models. Since we are interested in studying the spatially localized aggregation behavior as discussed in Section \ref{sec:blowup}, it is reasonable to consider the system \eqref{KS_rho} and \eqref{KS_c} in a large domain $\Omega=[-L/2,L/2]^d$ and assume Dirichlet boundary condition for the particle density $\rho$ and Neumann boundary condition for the chemical concentration $c$.
	
	As a discrete algorithm, we assume that the temporal domain $[0,T]$ is partitioned into time steps by $\{t_n\}_{n=0:n_T}$, where $t_0=0$ and $t_{n_T}=T$. We approximate the density $\rho$ using particles, i.e. 
	\begin{align}
		\rho_t \approx {\frac {M_0}{P}}\, \sum_{j=1}^{P} \delta (x - X^{p}_{t}),  \; \; P\gg 1,
		\label{rho_particle}
	\end{align}
	where $M_0$ is the conserved total mass (integral of $\rho$).
	For the chemical concentration $c$, we approximate it using a set of Fourier basis. \rev{In the following derivation, we assume $d=3$ for brevity, while the algorithms proposed work under any dimension.} Specifically, $c(\textbf{x},t)$ can be represented as a series expansion: 
	\begin{equation}
		\sum_{j,m,l \in \mathcal{H}}\, \alpha_{t;j,m,l}\, \exp(i 2\pi j \, x_1/L) \exp(i 2\pi m\, x_2/L)\exp(i2\pi l\, x_3/L), \label{trig_ser}
	\end{equation}
	where $\mathcal{H}$ denotes index set 
	\begin{align}\label{def:H}
		\mathcal{H}:=\{(j,m,l)\in \mathbb{N}^3: |j|,|m|,|l|\leq \frac{H}{2}\},
	\end{align}
	and $i=\sqrt{-1}$. 
	
	Then at $t_0=0$, we generate $P$ empirical samples $\{ X^{p}_{0}\}_{p=1:P}$ according to the initial condition of $\rho_0$, and set up $\alpha_{0;j,m,l}$ using the Fourier series representation of $c_0$. 
	
	For ease of presenting our algorithm, we will use a slight abuse of notation. We will 
	represent the density $\rho$ at time $t_n$ as $\rho_n={\frac {M_0}{P}}\, \sum_{p=1}^{P} \delta (x - X^{p}_{n})$, and the chemical concentration $c$ at time $t_n$ as 
	\[ c_n=\sum_{j,m,l \in \mathcal{H}}\, \alpha_{n;j,m,l}\, \exp(i 2\pi j \, x_1/L) \exp(i 2\pi m\, x_2/L)\exp(i2\pi l\, x_3/L).\]
  
	To discretize the time-stepping system (\ref{ppKS}) from $t_n$ to $t_{n+1}$, with $\rho_n$ and $c_{n-1}$ known, our algorithm, inspired by the operator splitting technique, consists of two sub-steps: updating the chemical concentration $c$ and updating the organism density $\rho$.

	\paragraph{Updating chemical concentration $c$}
	
	Let $\delta t=t_{n+1}-t_n>0$ be the time step. We discretize the equation for $c$ in \eqref{ppKS} in time using an implicit Euler scheme:
	\begin{align}
		\epsilon \, ( c_n - c_{n-1} )/\delta t
		= (\Delta - k^2 )\, c_n + \rho_n. \label{Euler4c} 
	\end{align}
	From \eqref{Euler4c}, we obtain the explicit formula for $c_n$ as:
	\begin{align}
		(\Delta - k^2 -\epsilon/\delta t )\, c_n
		= - \epsilon \, c_{n-1}/\delta t - \rho_n.  \label{Euler4c2}
	\end{align}
	It follows that:
	\begin{align}
		c_n=c(\textbf{x},t_n) &= - \mathcal{K}_{\epsilon,\delta t} \ast (\epsilon \, c_{n-1}/\delta t + \rho_n) 
		= - \mathcal{K}_{\epsilon,\delta t} \ast (\epsilon \, c(\textbf{x},t_{n-1})/\delta t +\rho(x,t_n))
		\label{IPS_Cstep}
	\end{align}
	where $\ast$ is the spatial convolution operator, and  $\mathcal{K}_{\epsilon,\delta t}$ is the Green's function of the operator $\Delta - k^2 -\epsilon/\delta t$. In the case of $\mathbb{R}^3$, the Green's function $\mathcal{K}_{\epsilon,\delta t}$ is given by:
	\begin{align}
		\mathcal{K}_{\epsilon,\delta t}
		=  \mathcal{K}_{\epsilon,\delta t}(\textbf{x})= -  
		\frac{\exp\{ - \beta |\textbf{x} | \}}{ 4 \pi |\textbf{x}|}. \quad
		\beta^2 = k^2 + \epsilon/\delta t. \label{GreenFunction3D}
	\end{align}
	The Green's function admits a closed-form Fourier transform,
	\begin{align}
		\mathcal{F}\mathcal{K}_{\epsilon,\delta t} (\omega)=-\frac{1}{|\omega|^2+\beta^2}.\label{GreenFunction3D_fourier}
	\end{align}
	
	For the term $-\mathcal{K}_{\epsilon,\delta t} \ast  c_{n-1}$ in \eqref{IPS_Cstep}, using Eq.\eqref{GreenFunction3D_fourier}, it is equivalent to modify the Fourier coefficients $\alpha_{j,m,l}$ to $\alpha_{j,m,l}/(4\pi^2j^2/L^2+4\pi^2m^2/L^2+4\pi^2l^2/L^2+\beta^2)$.  
	\medskip
	
	For the second term $\mathcal{K}_{\epsilon,\delta t} \ast \rho$, we first approximate $\mathcal{K}_{\epsilon,\delta t}$ with a $\cos$ series expansion. Then, according to the particle representation of $\rho$  in \eqref{rho_particle}, we have
    $$
    \begin{array}{ll}
		& (\mathcal{K}_{\epsilon,\delta t} \ast \rho)_{j,m,l} 
        \approx \nonumber \\
        & \frac{M_0}{P}\sum_{p=1}^P \frac{\exp(-\rev{i}2\pi j X^p_{n,1}/L-\rev{i}2\pi m X^p_{n,2}/L-\rev{i}2\pi l X^p_{n,l}/L)(-1)^{j+m+l}}{4\pi^2j^2/L^2+4\pi^2m^2/L^2+4\pi^2l^2/L^2+\beta^2}.
        \label{integralKrho}
         \end{array}
    $$
	Finally, we summarize the one-step update of the Fourier coefficients of 
    the concentration field  
    $c$ in  Algorithm \ref{alg:SIPM-C}.
	\begin{algorithm}[htbp]
		\KwData{Distribution $\rho_n$ represented by empirical samples $X_n$ , initial concentration $c_{n-1}$ represented by Fourier coefficients $\alpha_{n-1}$; \\}
		\SetAlgoLined
		\For{$(j,m,l)\in \mathcal{H}$}{
			$\alpha_{n;j,m,l}\gets \frac{\epsilon \alpha_{n-1;j,m,l}}{\delta t(4\pi^2j^2/L^2+4\pi^2m^2/L^2+4\pi^2l^2/L^2+\beta^2)}$\\
			$F_{j,m,l}\gets 0$.\\
			\For{$p=1$ \KwTo P}{
				$F_{j,m,l}\gets F_{j,m,l}+\exp(-\rev{i}2\pi j X^p_{n;1}/L-\rev{i}2\pi m X^p_{n;2}/L-\rev{i}2\pi l X^p_{n;3}/L)
				$
			}
			$F_{j,m,l}\gets F_{j,m,l}\frac{(-1)^{j+m+l}}{4\pi^2j^2/L^2+4\pi^2m^2/L^2+4\pi^2l^2/L^2+\beta^2}*\frac{M}{P}$
		}
		$\alpha_n\gets \alpha_n-F$\\
		\KwResult{Updated chemical concentration {field from input $c_{n-1}$ to $c_n$} via $\alpha_n$.}
		\caption{One step update of chemical concentration in SIPF}\label{alg:SIPM-C}
	\end{algorithm}
    \rev{\begin{remark}Note that \eqref{Euler4c}-\eqref{IPS_Cstep} are algebraic derivations. In the vanishing $\epsilon$ regime, the related semi-discrete systems tend to an elliptic equation. In terms of numerical implementation, the updates of Fourier coefficients with Green's kernel by \eqref{GreenFunction3D} and \eqref{integralKrho} are also consistent. Hence our method is robust in the vanishing $\epsilon$ regime.
    \end{remark}}
	
	\paragraph{Updating density of active particles $\rho$}
	In the one-step update of the density $\rho_n$ represented by particles $\{X_n^p\}_{p=1:P}$, we apply the Euler-Maruyama scheme to solve the SDE 
	\eqref{IPS_SDE}:
	\begin{align}
		X^{p}_{n+1} = X^{p}_{n} 
		+\chi \nabla_\textbf{x} c(X^{p}_{n},t_n)\delta t
		+ \sqrt{2\,\mu\, \delta t} \, N^p_n, 
		\label{OnestepX_EM}
	\end{align}
	where $N^p_n$'s are i.i.d. standard normal distributions with respect to the Brownian paths in the SDE formulation \eqref{IPS_SDE}. For $n>1$, substituting \eqref{IPS_Cstep} in \eqref{OnestepX_EM} gives:
	\begin{align}
		X^{p}_{n+1} = X^{p}_{n} 
		-\chi \nabla_\textbf{x} \mathcal{K}_{\epsilon,\delta t} \ast (\epsilon\, c_{n-1}(\textbf{x})/\delta t + \rho_n(\textbf{x}))|_{\textbf{x}=X^{p}_{n}}\delta t
		+ \sqrt{2\,\mu\, \delta t} \,  N^p_n, 
		\label{IPS_Xstep_full}
	\end{align}
	from which $\rho_{n+1}(\textbf{x})$ is constructed via (\ref{rho_particle}). 
	
	In such a particle formulation, the computation of spacial convolution is slightly different from the one in the update of $c$, namely \eqref{IPS_Cstep}. 
	
	For $\nabla_{\textbf{x}}\mathcal{K}_{\epsilon,\delta t} \ast \, c_{n-1}(X^p_n)$, to avoid the singular points of $\nabla_{\textbf{x}}\mathcal{K}_{\epsilon,\delta t}$, we evaluate the integral with quadrature points that are away from $0$. Precisely, we denote the standard quadrature points in $\Omega$ as 
	\begin{equation}
		x_{j,m,l}=( j\, L/H,m\, L/H, \revv{l}\, L/H), \label{grid}
	\end{equation}
	where $j$, $m$, $l$ are integers ranging from $-H/2$ to $H/2-1$. When computing the integral $\nabla_{\textbf{x}}\mathcal{K}_{\epsilon,\delta t} \ast \, c_{n-1}(X^p_n)$, we evaluate $\nabla_{\textbf{x}}\mathcal{K}_{\epsilon,\delta t}$ at $\{X^p_n+\bar{X}^p_n-x_{j,m,l}\}_{j,m,l}$ where a small spatial shift is given by 
 \begin{align}
\label{eqn:sp_shft}\bar{X}^p_n=\revv{\frac{L}{2H}+\lfloor \frac{X^p_n}{L/H}\rfloor \frac{L}{H}}-X^p,
 \end{align} and $c$ at $\{x_{j,m,l}-\bar{X}^p_n\}_{j,m,l}$ correspondingly. The latter one is computed by the inverse Fourier transform of the shifted coefficients, with $\alpha_{j,m,l}$ modified to $\alpha_{j,m,l}\exp(-i2\pi j\bar{X}^p_{n;1}/L- i2\pi m\bar{X}^p_{n;2}/L-i2\pi l\bar{X}^p_{n;3}/L) $, where $(\bar{X}^p_{n;i})$ denotes the $i$-th component of $\bar{X}^p_{n}$.
	
	The computation of the term  $\nabla_{\textbf{x}}\mathcal{K}_{\epsilon,\delta t} \ast \, \rho(X^p_n,t_{n})$
	is straightforward thanks to the particle representation of  $\rho(X^p_n,t_{n})$ in \eqref{rho_particle}: 
	\begin{align}
		\nabla_{\textbf{x}}\mathcal{K}_{\epsilon,\delta t} \ast \, \rho_n(X^p_n)=\int \nabla_{\textbf{x}}\, \mathcal{K}_{\epsilon,\delta t}(X^p_n-y)\rho(y)\approx \sum_{q=1,q\not=p}^P\frac{M}{P}\nabla_{\textbf{x}}\, \mathcal{K}_{\epsilon, \delta t}(X^p_n-X^q_n).
	\end{align}
	
	We summarize the one-step update (for $n>1$) of density in SIPF as in Algorithm \ref{alg:SIPM-rho}.
	
	\begin{algorithm}[htbp]
		\KwData{Distribution $\rho_n$ represented by empirical samples $X_n$, input: concentration $c_{n-1}$ represented by Fourier coefficients $\alpha_{n-1}$; \\}
		\SetAlgoLined
		
		\For{$p=1$ \KwTo P}{
			$X_{n+1}^p\gets X_{n+1}^p+\sqrt{2\mu\delta t}N$ where $N$ is a random generated standard normal distribution.\\
			\For{$q=1$ \KwTo P}{
				$X_{n+1}^p\gets X_{n+1}^p -\frac{\chi M \delta t}{P}\nabla_{\textbf{x}}\mathcal{K}_{\epsilon, \delta t}(X^p_n-X^q_n)$
			}
			$\bar{X}^p_n\gets\frac{H}{2L}+\lceil \frac{X^p_n}{H/L}\rceil \frac{H}{L}-X^p$\\
			
			\For{$(j,m,l)\in \mathcal{H}$}{
				$F_{j,m,l}\gets \nabla_{\textbf{x}}\mathcal{K}_{\epsilon,\delta t}(X^p_n+\bar{X}_n^p-x_{j,m,l}),\;\;$ {$x_{j,m,l}$ from Eq. (\ref{grid}) }\\
				$G_{j,m,l}\gets\alpha_{j,m,l}\exp(-i2\pi j\bar{X}^p_{n;1}/L- i2\pi m\bar{X}^p_{n;2}/L-i2\pi l\bar{X}^p_{n;3}/L)$}
			
			$\check{G}=iFFT(G)$\\
			
			$X_{n+1}^p\gets X_{n+1}^p-\epsilon \chi(F,\check{G}) \frac{L^3}{H^3} $, where $(\cdot,\cdot)\frac{L^3}{H^3}$ denote an inner product corresponding to $L^2(\Omega)$ quadrature.
		}
		\KwResult{{Output} $\rho_{n+1}$ represented by updated $X_{n+1}.$}
		\caption{One step update of density in SIPF}\label{alg:SIPM-rho}
	\end{algorithm}
	
	Combining \eqref{IPS_Cstep} and \eqref{IPS_Xstep_full}, we conclude that the recursion from 
    
    \noindent $(\{X^{p}_{n}\}_{p=1:P}, \rho_n(\textbf{x}),c_{n-1}(\textbf{x}))$ 
    to $(\{X^{p}_{n+1}\}_{p=1:P}, \rho_{n+1}(\textbf{x}),c_n(\textbf{x}))$ is complete. 
	We summarize the SIPF method in the following Algorithm \ref{alg:SIPM}.
	\begin{algorithm}[htbp]
		\KwData{Initial distribution $\rho_0$, initial concentration $c_0$; \\}
		\SetAlgoLined
		Generate $P$ i.i.d samples following distribution $\rho_0$, $X^1, X^2, \cdots X^P$.\\
		\For{$p\gets 1$ \KwTo $P$}{Compute $X^{p}_{1}$ by \eqref{OnestepX_EM}, with $c_{-1}=c_0$}
		
		Compute $c_1$ by Alg.\ref{alg:SIPM-C} with $c_0$ and $\rho_1=\sum_{p=1}^P\frac{M}{P}\delta_{X^p_1}$.\\
		\For{step $n\gets 2$
			\KwTo $N=T/\delta t$ }{{Compute $X_{n}$ by Alg.\ref{alg:SIPM-rho} with $\rho_{n-1}$ and $c_{n-2}$}
			
			Compute $c_{n}$ by Alg.\ref{alg:SIPM-C} with $c_{n-1}$ and $\rho_{n}=\sum_{p=1}^P\frac{M}{P}\delta_{X^p_n}$.}
		
		\caption{Stochastic Interacting Particle-{Field} Method}\label{alg:SIPM}
	\end{algorithm}

\rev{
\paragraph{Numerical Stability}
Under discretization $\mathcal{H}$, the update in Alg.\ref{alg:SIPM-C} is equivalent to,
\begin{align}
    \hat{c}_n(\omega)=\frac{1}{1+(|\omega|^2+k^2)\cdot \frac{\delta t}{\epsilon}}\hat{c}_{n-1}(\omega)+\frac{1}{|\omega|^2+k^2+\frac{\epsilon}{\delta t}}\hat{\rho}_n(\omega),
\end{align}
where for brevity we used the notation, $\omega=\frac{2\pi}{L}(n,m,l)$.
Given $|\hat\rho_n(\omega)|$ bounded by $M_\rho$,
\begin{align}
    |\hat{c}_n(\omega)| \leq M_c + A_c^n(|\hat{c}_0|-M_c),
\end{align}
where $A_c=\frac{1}{1+(|\omega|^2+k^2)\cdot \frac{\delta t}{\epsilon}}<1$ and $M_c = \frac{M_\rho}{|\omega|^2+k^2}$. Furthermore, in the discrete setting,
\begin{align}
    |\hat\rho_n(\omega)| = |\frac{M}{P}\sum_{p=1}^P \exp(-i2\pi j X^p_{n;1}/L-i2\pi m X^p_{n;2}/L-i2\pi l X^p_{n;3}/L)| \leq M.\label{est:rhon}
\end{align}
\eqref{est:rhon} implies that one can take $M_\rho=M$ (total mass) and validate that the stability of the update of $c$ in Alg.\ref{alg:SIPM-C} only requires $k>0$, independent of $\delta t$ and position of $X_n$.

In the one-step update of $\rho_n$ represented by $X_n$ via Alg.\ref{alg:SIPM-rho}, the numerical stability relates to the increment from drift terms, namely $\frac{\chi M \delta t}{P}\nabla_{\textbf{x}}\mathcal{K}_{\epsilon, \delta t}(X^p_n-X^q_n)$ and $\epsilon \chi(F,\check{G}) \frac{L^3}{H^3}$. The upper boundedness of the latter is ensured by the spatial shifting \eqref{eqn:sp_shft}. For the direct interaction term $\frac{\chi M \delta t}{P}\nabla_{\textbf{x}}\mathcal{K}_{\epsilon, \delta t}(X^p_n-X^q_n)$, due to the random normal update step, namely $X_{n+1}^p\gets X_{n+1}^p+\sqrt{2\mu\delta t}N$, the probability that $X^p_n$, $X^q_n$ coincides is zero. Practically, we clip $X_{n+1}$ to the range of domain $\Omega$.
}
 
	\rev{\begin{remark}
	    The convergence of the distribution of many-particle solution to the continuous PDE (the mean-field limit) in similar systems can be found in \cite{stevens2000derivation,chen2022mean}. In the mean-field regime, our discretization for the $c$ equation, namely \eqref{KS_c}, is implicit Euler in time and Fourier spectral discretization in space. We observe the first order convergence in time through numerical experiment. in Section \ref{sec:convergence-dt}. We will leave the full convergence analysis of the SIPF as a future work.
	\end{remark}}
\rev{\begin{remark}
     The memory usage throughout Alg.\ref{alg:SIPM} is $\mathcal{O}(P+\operatorname{Card} \mathcal{H})$, where $\operatorname{Card} \mathcal{H}$ denotes the total number of Fourier terms in the discretization \eqref{def:H}. The computational cost for each step of  Alg.\ref{alg:SIPM-rho}
and Alg.\ref{alg:SIPM-C} are $\mathcal{O}(P\operatorname{Card} \mathcal{H})$ and $\mathcal{O}(P^2 + P\operatorname{Card} \mathcal{H} \log(\operatorname{Card}\mathcal{H}))$. So the one step update in Alg.\ref{alg:SIPM} costs $\mathcal{O}(P^2 + P\operatorname{Card} \mathcal{H} \log(\operatorname{Card}\mathcal{H}))$. To be noted, the complete algorithm Alg.\ref{alg:SIPM} admits direct parallelization over particles in each step. \end{remark}}
	
	\section{Numerical Experiments}
	
	\subsection{Aggregation Behaviors}\label{sec:agrregation}
\begin{figure}[htbp]
		\centering
		\subfigure[$T=0$]{\includegraphics[width=0.32\linewidth]{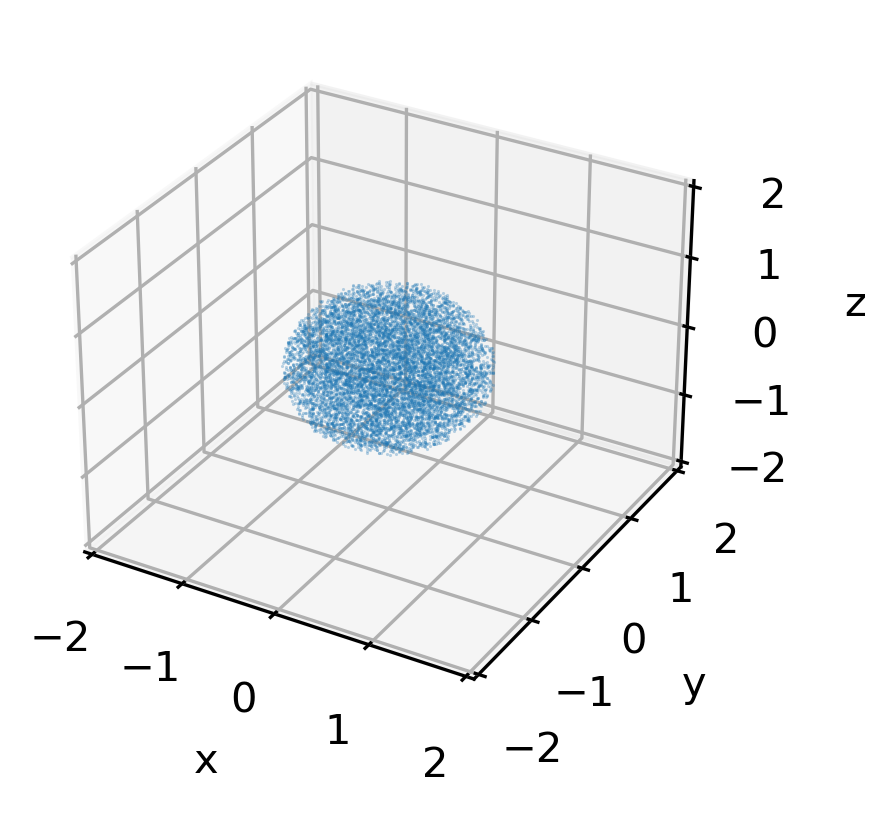}}
		\subfigure[$T=0.1$, $M_0=20$]{\includegraphics[width=0.32\linewidth]{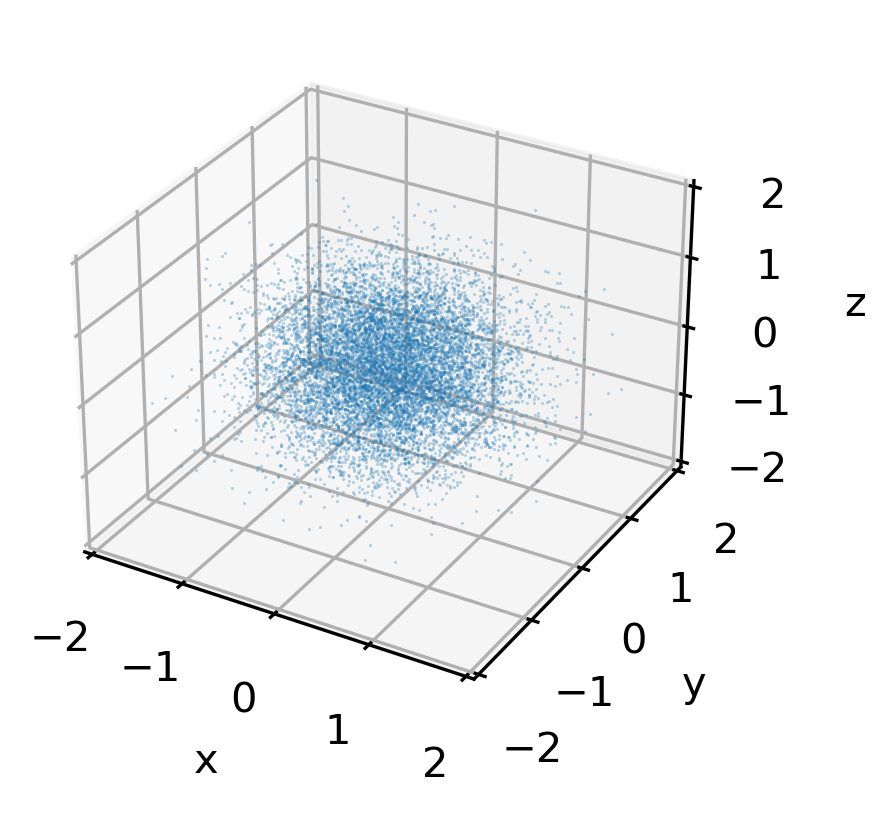}}
		\subfigure[$T=0.1$, $M_0=80$]{\includegraphics[width=0.32\linewidth]{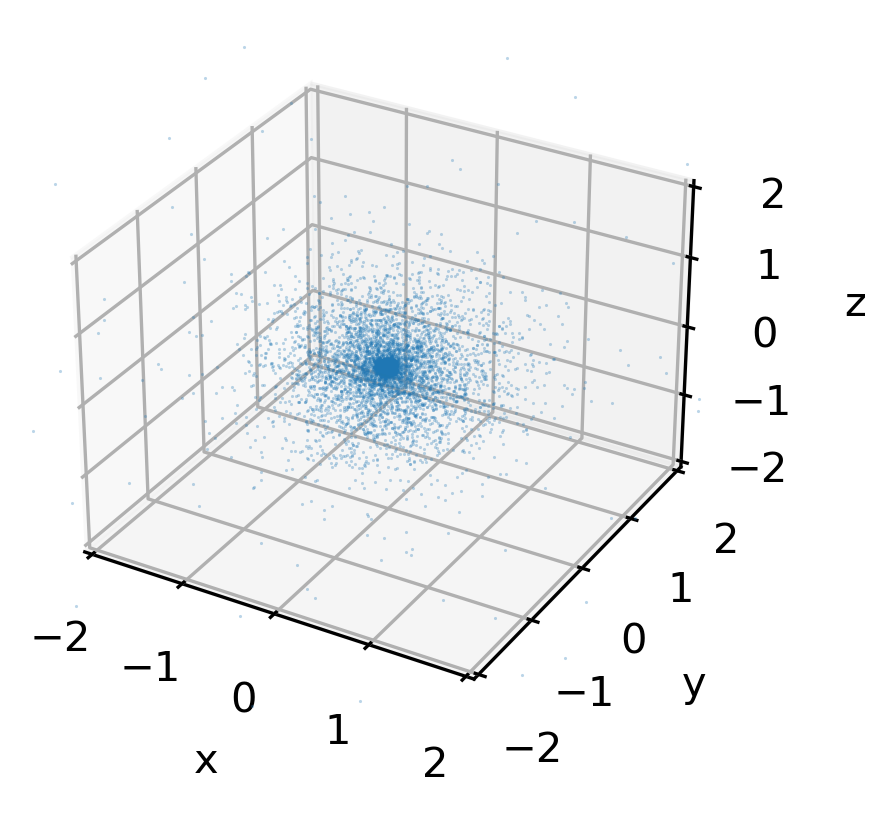}}
		
		\caption{Density $\rho$ approximated by empirical distribution at $T=0.1$: the mass effect on focusing.}
		\label{fig:eg1_rho}
	\end{figure}

	To illustrate the performance of the algorithm, we start with two examples. In both cases, the initial distribution $\rho_0$ is assumed to be a uniform distribution over a ball centered at $\begin{bmatrix} 0\\0\\0 \end{bmatrix}$ with radius $1$, as shown in Fig. \ref{fig:eg1_rho}(a). Also, in both cases, we assume the following model parameters,
	\begin{align}\label{eg_para}
		\mu=\chi=1, \quad \epsilon=10^{-4} \text{ and } k=10^{-1}.
	\end{align}
	for the fully parabolic KS model \eqref{ppKS}. These parameter choices in \eqref{eg_para} are made so that the model exhibits comparable behavior to the corresponding parabolic-elliptic KS system whose blow-up behavior is known. For the first example, the total mass is chosen to be $M_0=20$, while for the second one, the mass is $M_0=80$. 
	
	In the numerical computation of both examples, we use $H=24$ Fourier basis in each spatial dimension to discretize the chemical concentration $c$, and use $P=10000$ particles to represent the approximated distribution $\rho$. The computational domain is $\Omega=[-L/2, L/2]^3$ with $L=8$. We compute the evolution of $c$ and $\rho$ using Algorithm \ref{alg:SIPM} with a time step of $\delta t=10^{-4}$ up to $T=0.1$. The total computation time is $284$ seconds on a workstation with one RTX2070 Super GPU card.

	In Fig. \ref{fig:eg1_rho}, we plot the distribution $\rho$ by empirical samples, at the starting time $T=0$ and final computation time $T=0.1$.  In Fig. \ref{fig:eg1_rho}(b), we observe the diffusive behaviors compared to the initial distribution shown in Fig. \ref{fig:eg1_rho}(a). While in Fig. \ref{fig:eg1_rho}(c), we increase the total mass from $M_0=20$ to $M_0=50$, we can see particles become concentrated at the origin, which indicates the possible blow-up of the continuous system.
	
	\begin{figure}[htbp]
		\centering
		\subfigure[$M_0=20$]{\includegraphics[width=0.45\linewidth]{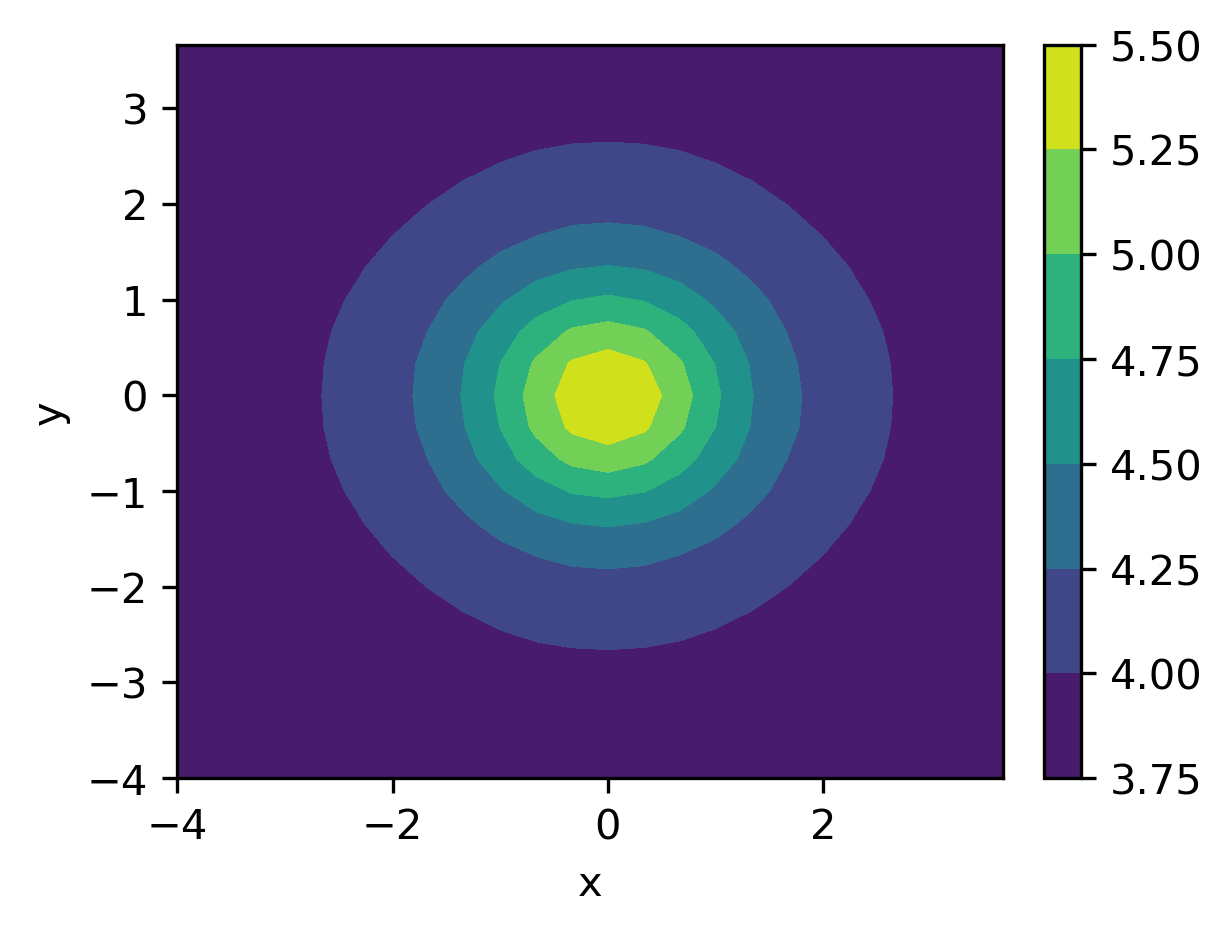}}
		\subfigure[$M_0=80$]{\includegraphics[width=0.45\linewidth]{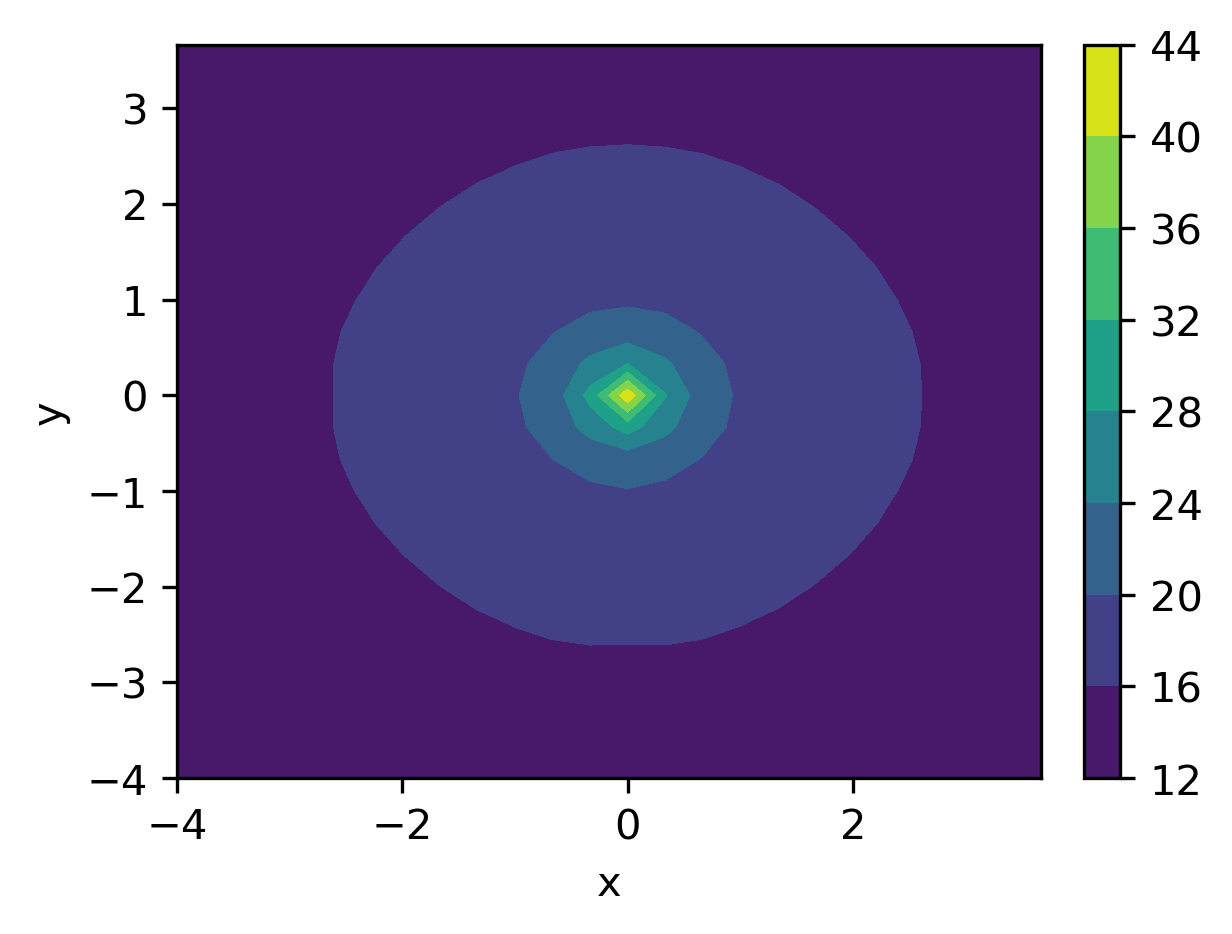}}
		
		\caption{Chemical concentration $c$ at final time $T=0.1$, sliced at $z=0$: \rev{the mass effect on focusing.}}
		\label{fig:eg1_c}
	\end{figure}
	
	In Fig. \ref{fig:eg1_c}, we present the chemical concentration $c$ at the final time $T=0.1$ and third component $z=0$ for various initial total masses $M_0$. By comparing the sub-figures, we can see that in the case of a large total mass, $c$ exhibits a sharp profile at the origin. This behavior, \rev{together with the near singular behavior of $\rho$, is shown in Fig.\ref{fig:eg1_rho} (b-c)}, 
    indicating a possible blow-up. 
	
	\begin{figure}[htbp]
		\centering
		\subfigure[$M_0=20$]{\includegraphics[width=0.45\linewidth]{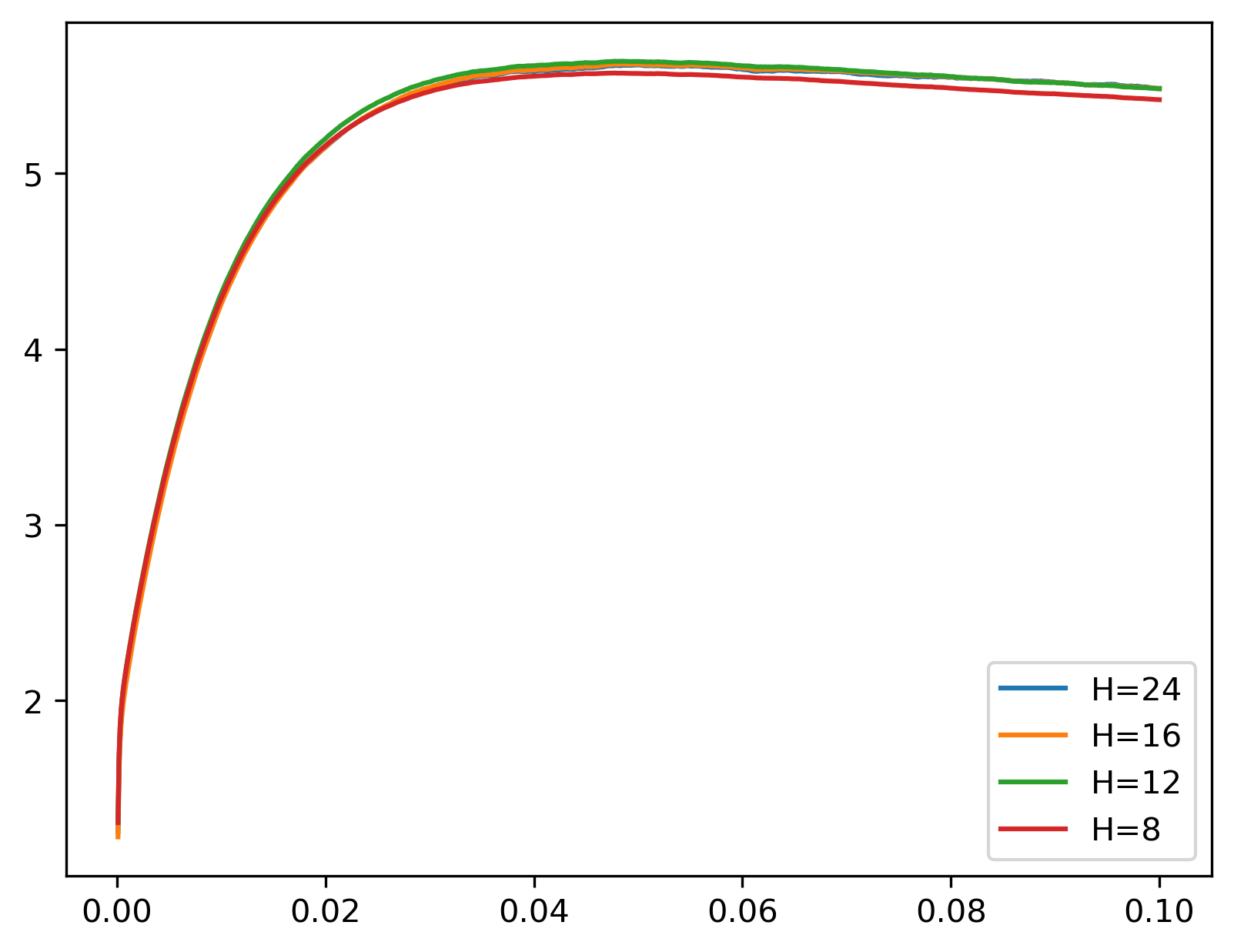}}
		\subfigure[$M_0=80$]{\includegraphics[width=0.45\linewidth]{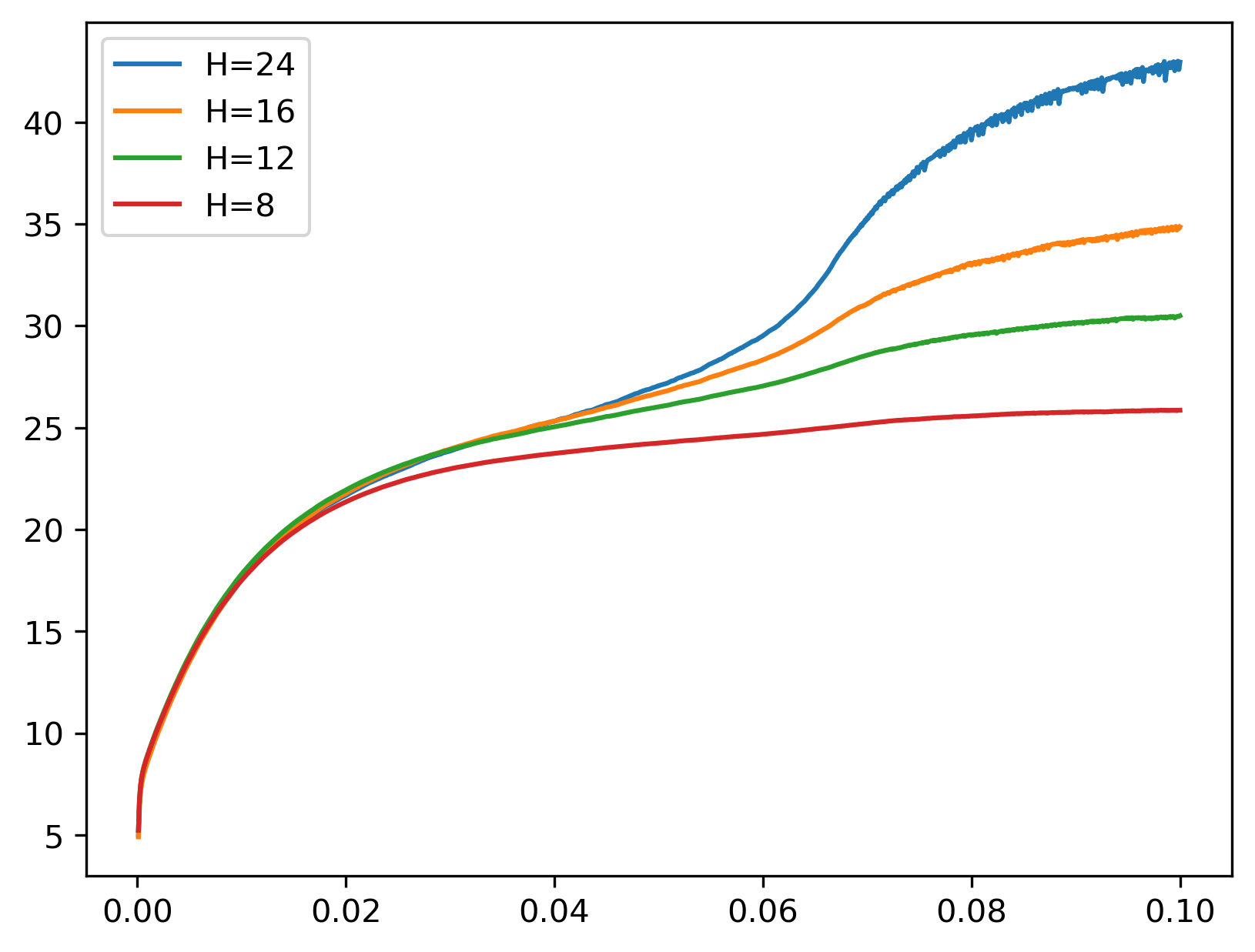}}
		
		\caption{Maximum of chemical concentration $c$ vs computation time $T$ with different total mass $M_0$: \rev{identifying blow up by refining the discretization. }}
		\label{fig:eg1_cmax}
	\end{figure}
	Furthermore, if we assume, there exists a self-similar profile of $\rho$ at origin as discussed in \cite{souplet2019blow} and Section \ref{sec:blowup}, namely, $\rho(x,t)\sim\frac{1}{|x|^2}$,
	the Fourier coefficients of the chemical concentration $c$ have the following asymptotics according to equation \eqref{ppKS}: 
	\begin{align}
		\mathcal{F}c(\omega)\sim \frac{1}{|\omega|^2+k^2} \hat{\rho}\sim\frac{1}{(|\omega|^2+k^2)|\omega|}.
	\end{align}
	Then, the maximum of $c$ in the computation shall vary with the discretization parameter $H$. Specifically, we observe that at the origin, 
	\begin{align}
		c(0)\sim \int \frac{1}{(|\omega|^2+k^2)|\omega|}e^{i\omega x}d\omega |_{x=0} = \int \frac{1}{(|\omega|^2+k^2)|\omega|}d\omega. \label{eq:asymptC0}
	\end{align}
	In practical discretization, the range of the integral \eqref{eq:asymptC0} is determined by the maximum frequency, which can be expressed as  $[-\frac{\pi}{L}(\frac{H}{2}-1),\frac{\pi}{L}\cdot\frac{H}{2}]^3$. Then, for the type of
	blow-up with a profile proportional to $\frac{1}{|x|^2}$, 
	\begin{align} \label{eq:asymptC1}
		\|c\|_{\infty}=\mathcal{O}(\ln(H)).
	\end{align}
	A similar analysis shows that for the type of blow-up with a profile proportional to $\delta(x)$, we have the following asymptotic relation:
	\begin{align}\label{eq:asymptC2}
		\|c\|_{\infty}=\mathcal{O}(H).
	\end{align}
	In Fig. \ref{fig:eg1_cmax}, we show the maximum value of $c$ as a function of computational time $T$ for different numbers of Fourier modes $H$ and total mass $M_0$. In the case of a possible blow-up (Fig. \ref{fig:eg1_cmax}(b)), we can see that the maximum of $c$ varies dramatically for different values of $H$. This variation can be used as an indicator of a possible blow-up, which will be further investigated in the following experiments.
	 
\begin{figure}[htbp]
		\centering
		\includegraphics[width=0.45\linewidth]{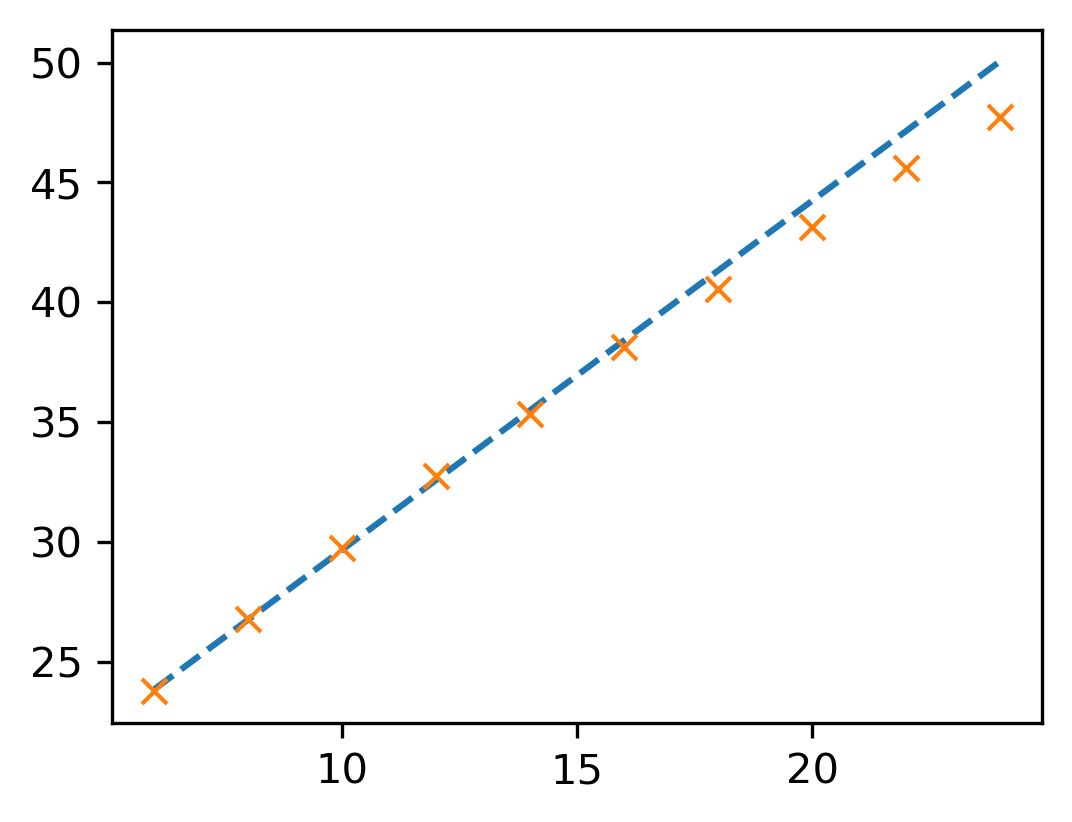}
		\caption{\rev{Maximum of $c$ scales linearly with the number of Fourier modes $H$ (in each dimension) under total mass $M_0=80$ (super critical): a $\delta$ type blow up.}}
		\label{fig:eg1_cmax_vs_H}
	\end{figure}

	
	Additionally, in the case where $M_0=80$ and $T=1$ are chosen to achieve numerical stability for $\|c\|_\infty$, we conducted tests for $H$ ranging from $8$ to $24$. In Fig. \ref{fig:eg1_cmax_vs_H}, we plot $\|c\|_\infty$ against $H$ and observe that the maximum value of $c$ grows nearly linearly with respect to $H$. This further supports our previous observation that the maximum of $c$  depends on the discretization parameter $H$.
	
	\begin{remark}
		Similar ideas for detecting blow-ups by comparing maximum values computed under different discretizations can be found in the literature on the finite volume approach to 2D KS systems. For example, in \cite{chertock2018high}, the $\delta$-type singularities in the 2D system are identified when $\|\rho\|_\infty=\mathcal{O}(\frac{1}{\Delta x \Delta y})$. 
	\end{remark}
	\subsection{Numerical Convergence}\label{sec:convergence-dt}
	\rev{In this subsection, we validate the convergence of the algorithms numerically. Throughout this subsection, we consider }  the same initial condition ($\rho$ and $c$ at $t=0$) and physical parameters as in the first example. \rev{For the parameters in the discretization, we will take uniform time step $\delta_{ref}=2^{-11}T$}, the number of Fourier modes in each dimension as $H=24$, the number of particles as $P=10000$, and the computational domain as $\Omega=[-L/2,L/2]^3$ with $L=8$. 
	Additionally, we set $M_0=80$ and $T=0.01$ when the system has not formed any singularities (as shown in Fig.  \ref{fig:eg1_cmax}(b)).

	For the investigation for convergence of $\delta t$, we consider the time step $\delta t$ in the range between $2^{-8}T$ to $2^{-4}T$, and we take $\delta t_{ref}=2^{-11}T$ as the reference solution. 
	
	Another identity that can be used to validate the accuracy of the computation for the system \eqref{KS_rho}-\eqref{KS_c} is the total concentration of the chemical attractant $c$ at a given time $t$, which is given by 
	\begin{align}\label{eq:cint}
		c_{\bf 0}(t):=\int_\Omega c(x,t)dx.
	\end{align}
	By integrating both sides of \eqref{KS_c} over physical domain $\Omega$, we have,
	\begin{align}
		\epsilon \frac{d}{dt}c_{\bf 0}(t)=-k^2c_{\bf 0}(t) + M_0,
	\end{align}
	which implies, 
	\begin{align}\label{eq:cint_groundtruth}
		c_{\bf 0}(t)=(c_{\bf 0}(0)-\frac{M_0}{k^2})e^{-\frac{k^2}{\epsilon}t}+\frac{M_0}{k^2}.
	\end{align}

	In Fig. \ref{fig:eg_conv_dt}, we present both types of error to verify the convergence of the algorithm. The first type of error is the $L_2$ relative error of the chemical concentration $c$ at the final time $T$ compared to the reference solution. The second type of error is the mean squared relative error of the total concentration \eqref{eq:cint} over the interval $[0,T]$ compared to the ground truth \eqref{eq:cint_groundtruth}. To calculate the mean squared relative error of the total concentration $c_{\bf 0}$, we use the following formula:
	\begin{align}
		c_{\bf 0}\ \text{error} = \sqrt{
			\frac{1}{T}\int_0^T \Big(\frac{|\tilde{c}_{\bf 0}(t)-c_{\bf 0}(t)|}{|c_{\bf 0}(t)|}\Big)^2 dt}\label{def:c0_err}
	\end{align}
	where $\tilde{c}_{\bf 0}$ is the approximated value of ${c}_{\bf 0}$ obtained using the SIPF algorithm with a time step $\delta t$. In addition, we fit the slope of error versus $\delta t$ in the logarithmic scale and find that the $L_2$ relative error $e(\delta t)$ follows an approximate first-order convergence rate, with $e(\delta t)=\mathcal{O}(\delta t^{1.011})$. 
	This indicates that the algorithm is approximately first order in time.
    
 \rev{
	\medskip
	\begin{minipage}[b]{0.45\linewidth}
	\centering
		\includegraphics[width=0.8\textwidth]{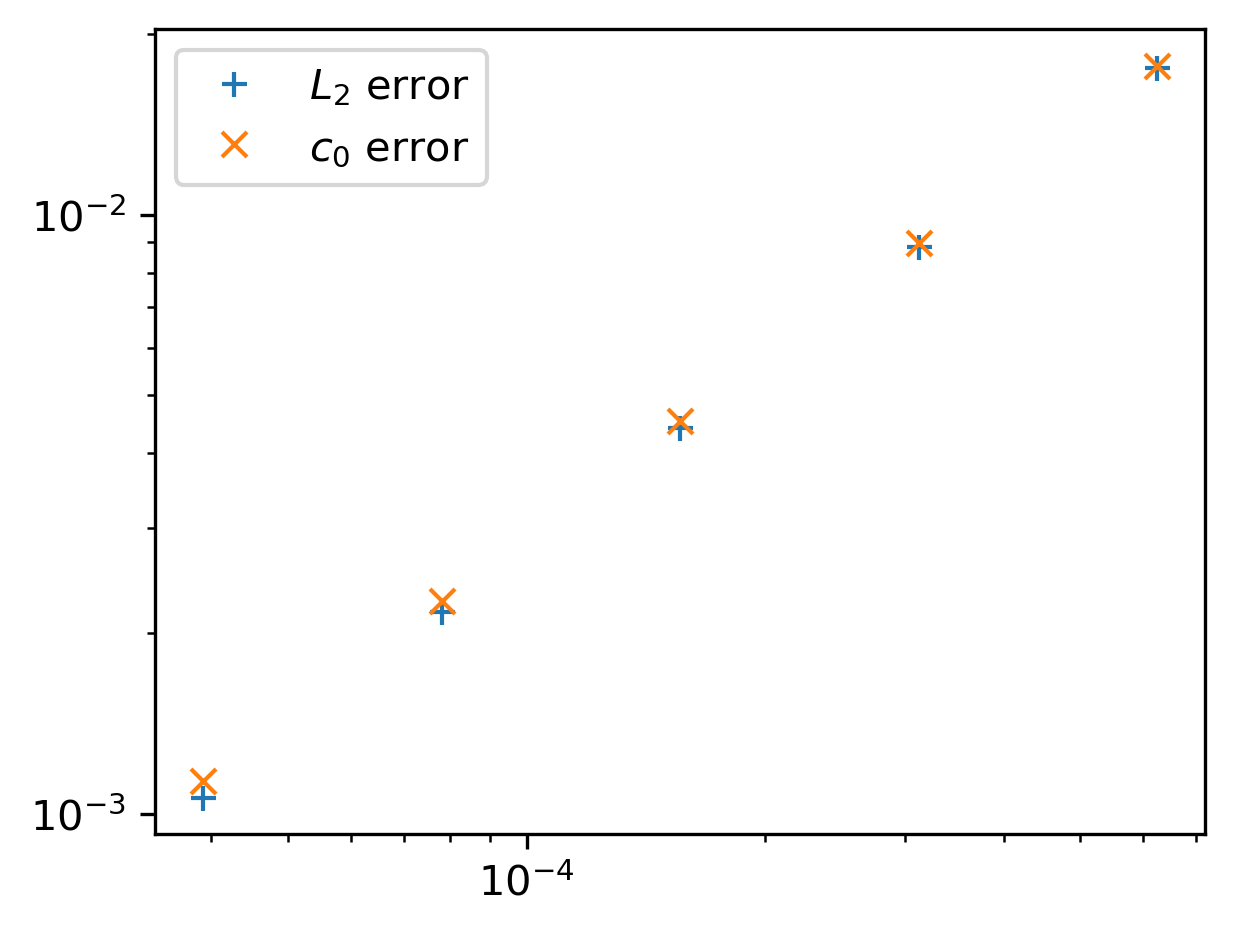}
		\captionof{figure}{Relative errors of $c$ vs. $\delta t$, compared with $\delta= 2^{-11}\times 0.01$. \rev{Fitted rate: $e(\delta t)=\mathcal{O}(\delta t^{1.011})$}.}
		\label{fig:eg_conv_dt}
 \end{minipage}\hfill
\begin{minipage}[b]{0.45\linewidth}
\centering
     \begin{tabular}{ccc}
     \toprule
         $\delta t$ & $L_2$ error & $c_{\bf 0}$ error \\
         \midrule
        $2^{-4}\times 0.01$  & $0.01756$  & $0.01773$\\
        $2^{-5}\times 0.01$  &  $0.00884$& $0.00898$\\
         $2^{-6}\times 0.01$ & $0.00440$ &$0.00452$\\
         $2^{-7}\times 0.01$  & $0.00218$  & $0.00227$\\
         $2^{-8}\times 0.01$  & $0.00106$ & $0.00113$\\
         \bottomrule
     \end{tabular}
     \captionof{table}{Relative errors of $c$ vs. $\delta t$, from data of Fig.\ref{fig:eg_conv_dt}.}
     \label{tab:my_label}
	 \end{minipage}

    In addition, we investigate the convergence of the algorithm with respect to other parameters of discretization. To this end, we keep the reference setting as Fig. \ref{fig:eg_conv_dt}, and alter Fourier mode number $H$ or particle number $P$ to compare with the reference solution. In Fig. \ref{fig:eg_conv_J} and Tab.  \ref{tab:err_J} we present the $L_2$ relative error of $c(\cdot,T)$ with varying $P=100,\ 200,\ 400,\ 800,\ 1600$ comparing with reference solution by $P=10000$. The fitted convergence rate is $e(P)=\mathcal{O}(P^{-0.287})$, significantly smaller than standard Monte Carlo. We conjecture that this is due to particle interaction and leave it for a future study. 
    
    In Fig. \ref{fig:eg_conv_H} and Tab. \ref{tab:err_H} we present the $L_2$ relative error of $c(\cdot,T)$ with varying $H=6,\ 8,\ 10,\ 12,$ comparing with reference solution by $H=24$. The fitted convergence rate is $e(H)=\mathcal{O}(H^{-0.659})$.
    
    We report that the convergent behavior of $c_{\bf 0}$  while varying $H$ or $P$ is similar. 

       \begin{minipage}[b]{0.45\linewidth}
	\centering
		\includegraphics[width=0.8\textwidth]{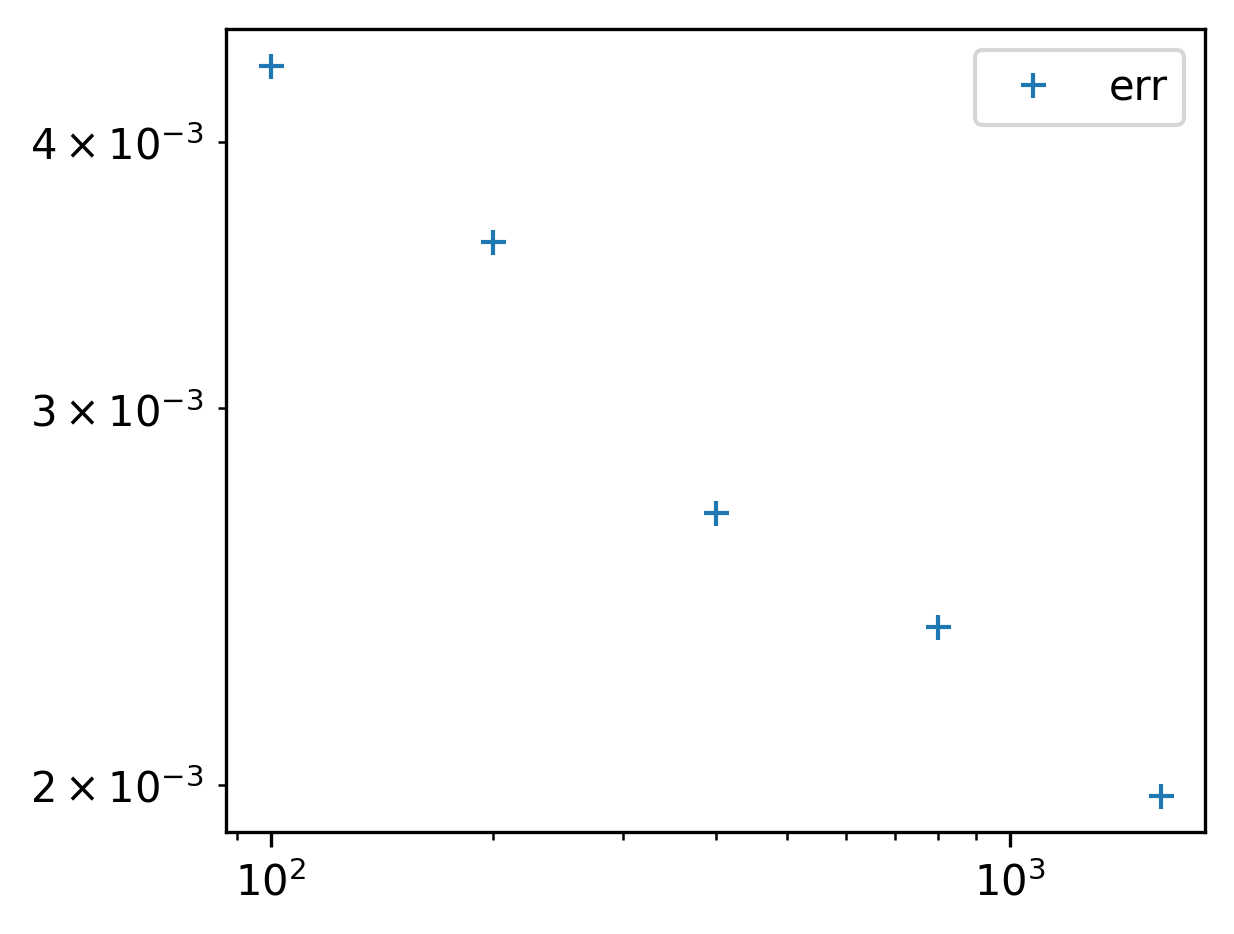}
		\captionof{figure}{Relative errors of $c$ vs. $P$, compared with $P=10000$. {.}}
		\label{fig:eg_conv_J}
 \end{minipage}\hfill
     \begin{minipage}[b]{0.45\linewidth}
\centering
     \begin{tabular}{cc}
     \toprule
         $P$ & $L_2$ error  \\
         \midrule
        $100$  & $0.00434$ \\
        $200$  &  $0.00359$\\
         $400$ & $0.00268$\\
         $800$  & $0.00237$ \\
         $1600$ & $0.00198$ \\
         \bottomrule
     \end{tabular}
     \captionof{table}{Relative errors of $c$ vs. $P$, from data of Fig. \ref{fig:eg_conv_J}.}
     \label{tab:err_J}
	 \end{minipage}

     \begin{minipage}[b]{0.45\linewidth}
	\centering
		\includegraphics[width=0.8\textwidth]{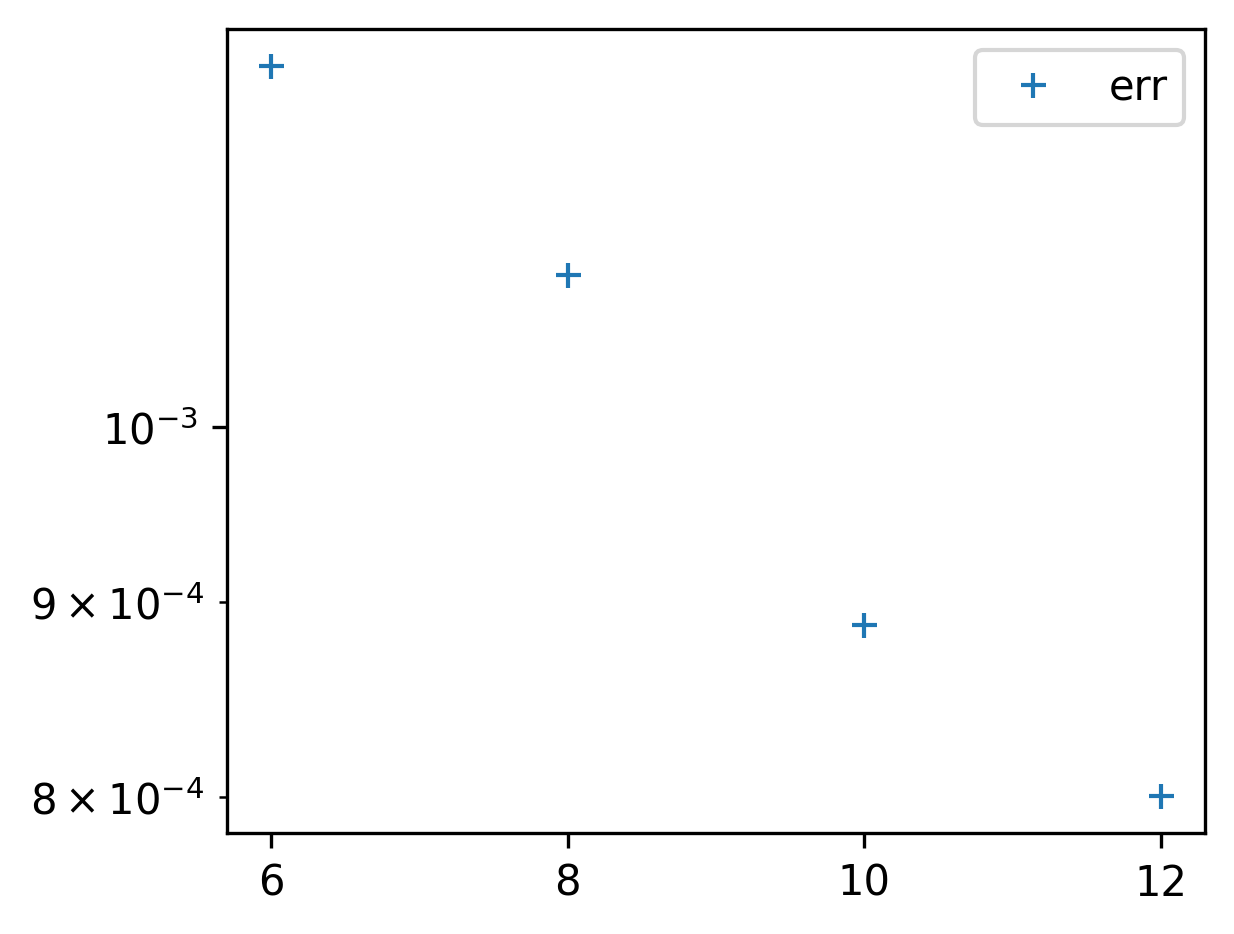}
		\captionof{figure}{Relative errors of  $c$ vs. $H$, compared with $H=24$. {Fitted rate: $e(H)=\mathcal{O}(H^{-0.659})$ }}
		\label{fig:eg_conv_H}
 \end{minipage}\hfill
     \begin{minipage}[b]{0.45\linewidth}
\centering
     \begin{tabular}{cc}
     \toprule
         $H$ & $L_2$ error  \\
         \midrule
        $6$  & $0.00124$ \\
        $8$  &  $0.00109$\\
         $10$ & $0.00089$\\
         $12$  & $0.00080$ \\
         \bottomrule
     \end{tabular}
     \captionof{table}{Relative errors of  $c$ vs. $H$, from data of Fig. \ref{fig:eg_conv_H}.}
     \label{tab:err_H}
	 \end{minipage}
}	
	\subsection{Blow up behaviors}\label{sec:blowup_num}
	As mentioned in Section \ref{sec:blowup}, it is well-known that $8\pi$ is the critical mass for the simplest two-dimensional parabolic-elliptic KS system \eqref{ppKS2Dsimple}.  Specifically, the system exhibits the following dichotomy:
	\begin{enumerate}
		\item If $M_0<8\pi$, the system has a global smooth solution.
		\item If $M_0>8\pi$, the system has no global smooth solutions and can exhibit blow-up behavior.
	\end{enumerate}
	
	In the case of a fully parabolic KS system or the specific parabolic-elliptic KS system \eqref{ppKS2Dsimple} with passive advection, there is no known variance identity similar to \eqref{simple_mass_result}. Therefore, to investigate the possible blow-up behaviors, numerical computation becomes necessary. By utilizing the asymptotics described in \eqref{eq:asymptC1} and \eqref{eq:asymptC2}, we can conduct tests for two scenarios: $H=24$ and $H=12$. In these examples, we will compare the $\|c\|_\infty$ values to detect any potential blow-up occurrences.  
	
\begin{figure}[htbp]
		\centering
		\subfigure[$\frac{|c|_{\infty, H=24}}{|c|_{\infty, H=12}}$ vs. computation time $T$.]{\includegraphics[width=0.32\linewidth]{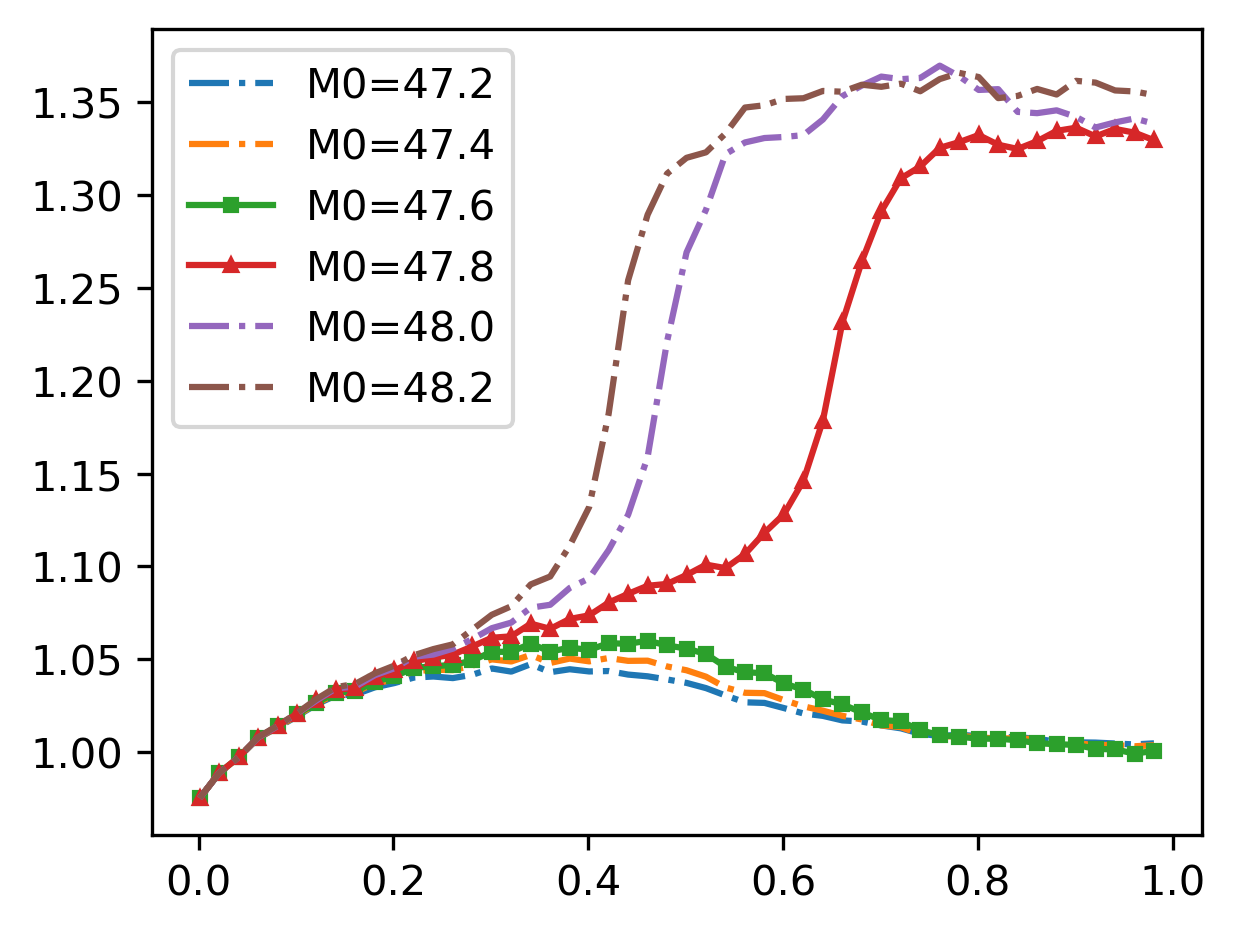}}
		\subfigure[$\frac{|c|_{\infty, H=24}}{|c|_{\infty, H=12}}$ at $T=1$ vs. $M_0$.]{\includegraphics[width=0.32\linewidth]{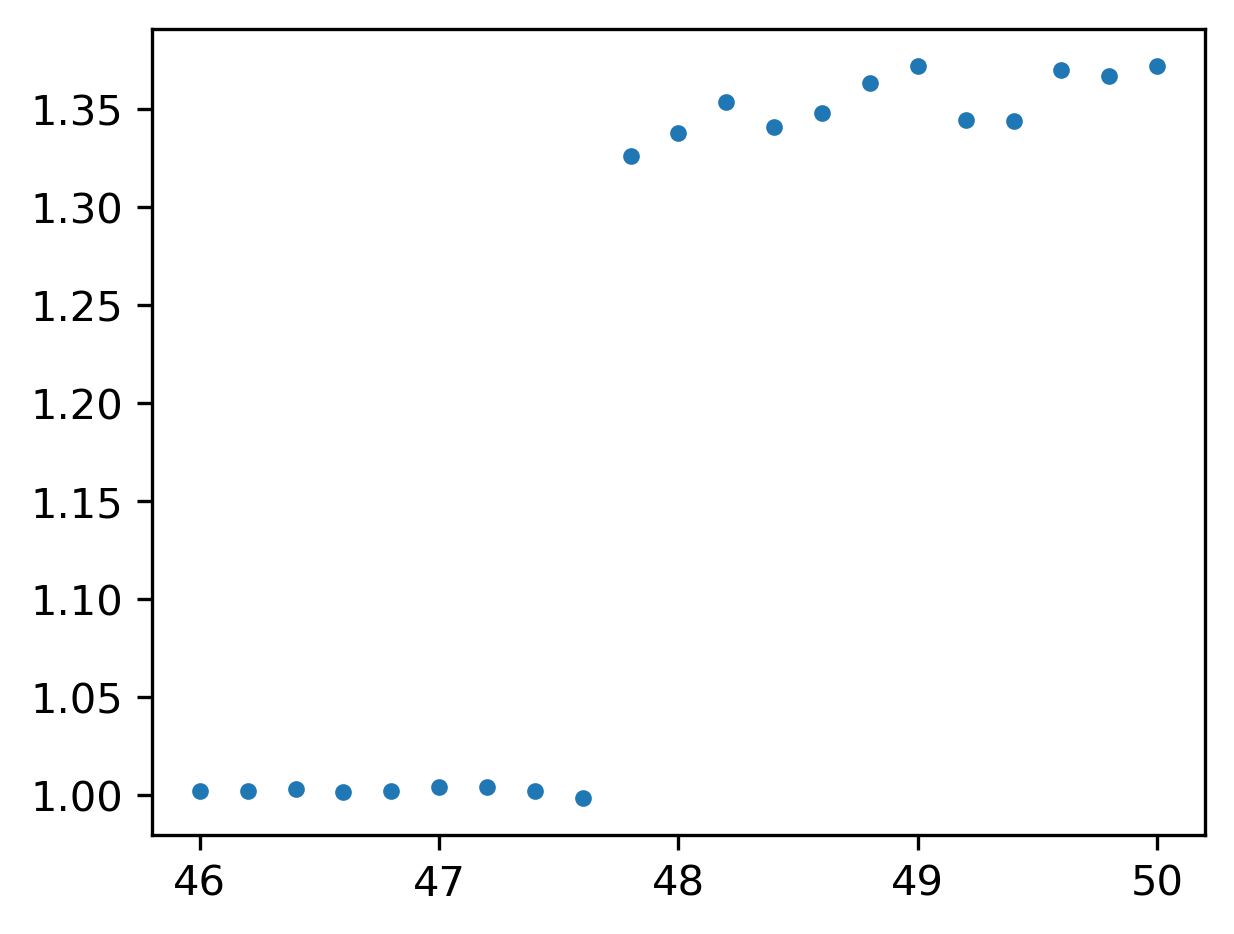}}
		\subfigure[$|c|_{\infty, FDM}$ at $T=1$ vs. $M_0$.  $\times$ denotes numerically unstable result.]{\includegraphics[width=0.32\linewidth]{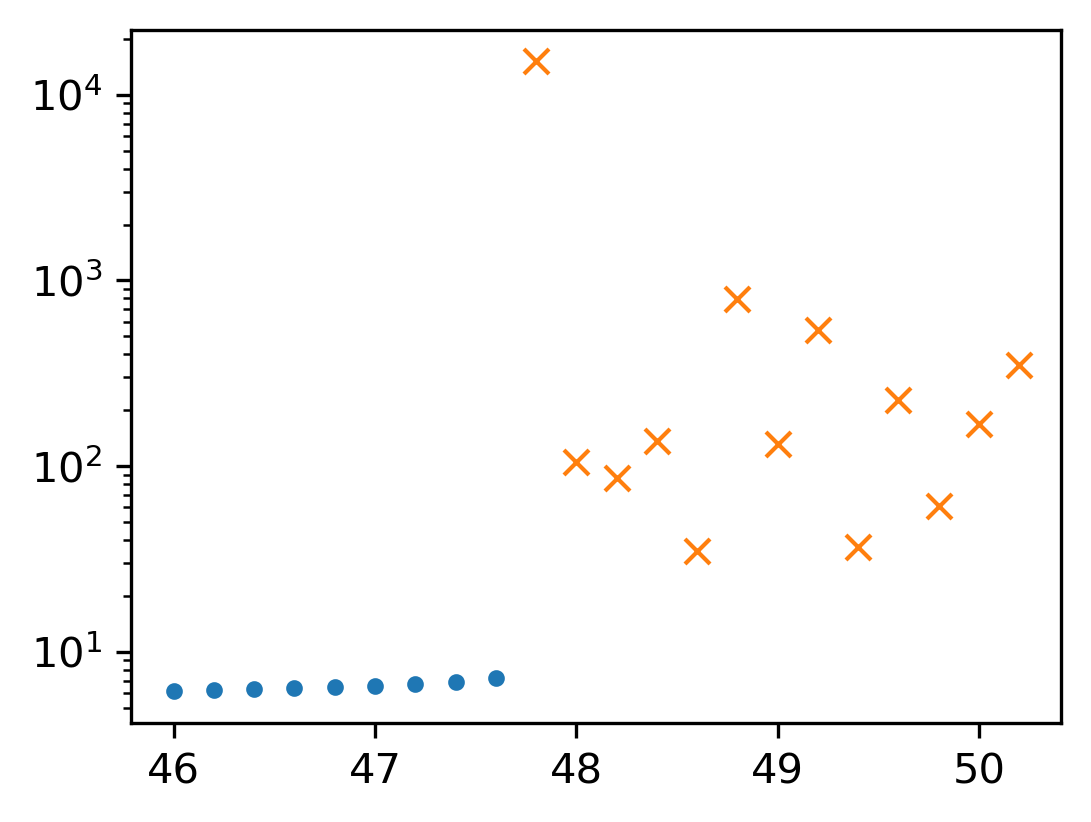}}
		\caption{Ratio of $|c|_\infty$'s from 2 SIPF runs with $H=24$ and $H=12$ and $|c|_\infty$ from FDM run: \rev{in radial case, the proposed SIPF can identify the same critical mass as FDM simulation.}}
		\label{fig:eg_mass}
	\end{figure}

	\paragraph{Mass dependence}
	We begin by investigating the critical mass $M_0$, which plays a central role in the dichotomy of the simple 2D parabolic elliptic system \eqref{ppKS2Dsimple}. To this end, we initialize the algorithm with a uniform distribution over the unit ball centered at the origin and $c(0,x)=0$. We then apply the algorithm with two different values of $H$ to compute the density and chemical concentration until $T=1$. To identify the possible blow-up, we compute the ratio of $|c|_\infty$ between the two cases. In Fig.\ref{fig:eg_mass}(a), we present the ratio, $\frac{|c|_{\infty, H=24}}{|c|_{\infty, H=12}}$, over time for various values of $M_0$. We observe a sharp increase in the ratio when a potential blow-up forms for $M_0\geq 47.6$. Fig.\ref{fig:eg_mass}(b) presents the ratio at the final time $T=1$, indicating that the critical mass of the aforementioned initial condition should fall between $47.6$ and $47.8$.
	
	In addition to the SIPF algorithm, we also present the numerical results obtained using the finite difference method (FDM). We note that the KS system \eqref{KS_rho}-\eqref{KS_c} admits radial solutions when given constant scalar physical parameters \eqref{eg_para} and a radially symmetric initial condition $(\rho_0,c_0)$. Therefore, we re-write the system in the radial coordinate,
	\begin{align}
		\rho_t&=\mu(\rho_{rr} + \frac{2}{r}\rho_r)-\chi(\rho_r c_r +\rho(c_{rr}+\frac{2}{r}c_r),\label{radial_rho}\\
		\epsilon c_t&=c_{rr}+\frac{2}{r}c_r-k^2c+\rho,
		\label{radial_c}\\
		r &\in \mathbb{R}^+, \quad t \in[0,T].
	\end{align}
	To formulate a finite difference scheme, we consider the system  \eqref{radial_rho}-\eqref{radial_c} on the domain $[0,20]$ with the Neumann boundary condition. We use a uniform partition with $N_r=2\times 10^5$ intervals for spatial discretization. For the temporal domain, we employ a backward Euler scheme with a time step of $\delta t=10^{-5}$. This discretization method requires a comparable computational time (approximately $150$ seconds) to the proposed SIPF. In Fig.\ref{fig:eg_mass}(c), we present the maximum value of $c$ over time for various initial mass $M_0$,  denoted as
	\begin{align}
		|c|_{\infty, FDM} = \sup_{t,r}|c(r,t)|.
	\end{align}
	We found that the FDM exhibits numerical instabilities for initial masses between $47.6$ and $47.8$, which matches the prediction made by the SIPF algorithm in Fig.\ref{fig:eg_mass}(b). This example further validates the accuracy of the proposed SIPF algorithm. It is worth noting that our SIPF algorithm applies readily to more general (non-radial) KS systems, whereas FDM in 3D with a fine uniform mesh will be computationally much more expensive.
	
\begin{figure}[htbp]
		\centering
		\subfigure[$\frac{|c|_{\infty, H=24}}{|c|_{\infty, H=12}}$ vs. computation time $T$.]{\includegraphics[width=0.45\linewidth]{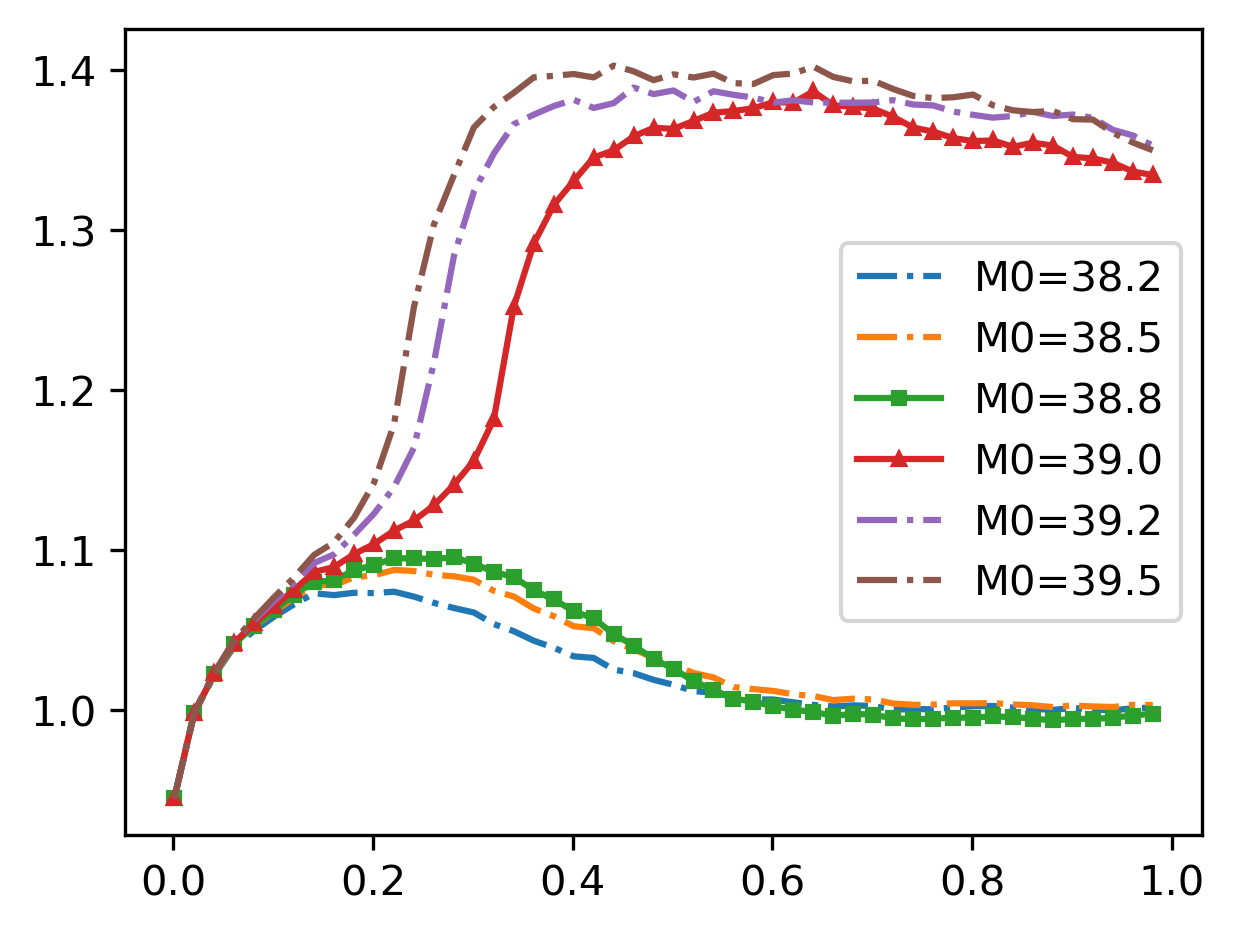}}
		\subfigure[$\frac{|c|_{\infty, H=24}}{|c|_{\infty, H=12}}$ at $T=1$ vs. $M_0$.]{\includegraphics[width=0.45\linewidth]{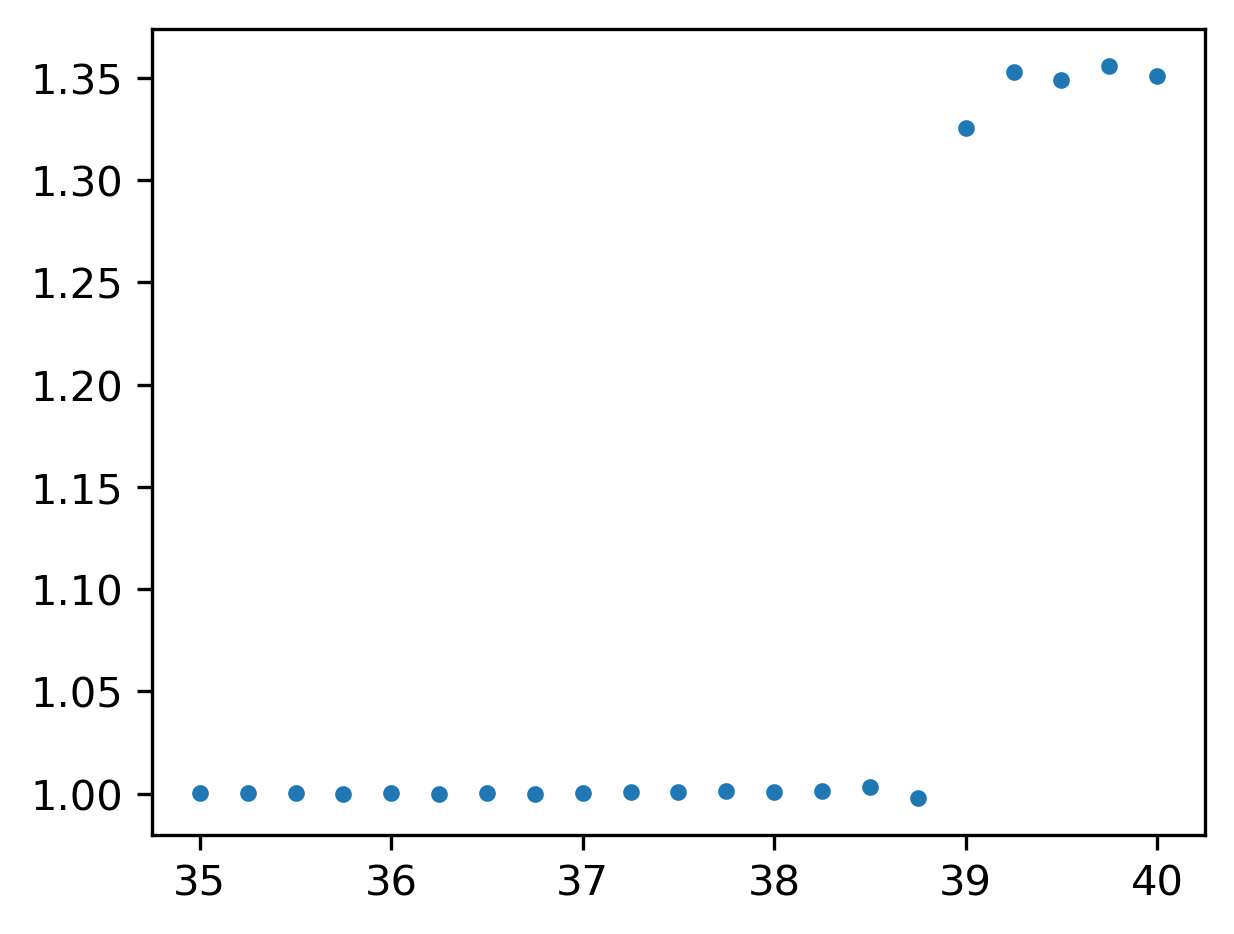}}
		
		\caption{Ratio of $|c|_\infty$'s from 2 runs with $H=24,12$; particles stay within initial radius $0.8$: \rev{the more concentrated initial, the smaller critical mass.}}
		\label{fig:eg_mass1}
	\end{figure}

	\paragraph{Geometry dependence}
	In contrast to the simplest parabolic-elliptic KS system \eqref{ppKS2Dsimple}, where the total mass is the sole determinant of the aggregation behavior, we have experimentally observed that the critical mass can vary for different initial distributions of $\rho$. For instance, in an experiment aimed at identifying the critical mass (as shown in Fig.\ref{fig:eg_mass}), we replaced the initial distribution with a uniform distribution on a ball centered at the origin with a radius of $\textbf{0.8}$. With a more concentrated initial distribution, we found that the critical mass of the system decreases. To be more precise, in Fig.\ref{fig:eg_mass1}(a), we present the ratio of $|c|_{\infty}$ for various total masses $M_0$ as a function of computational time $T$. We can see a significant change in the ratio when the total mass $M_0$ is large enough ($M_0\geq39$),  indicating the formation of potential singularities. Conversely, for relatively small values of $M_0$ ($M_0\leq 38.8$), the ratio remains stable around $1$ throughout the computational time. In Fig.\ref{fig:eg_mass1} (b), we present the ratio at the final time $T=0.1$ as a function of the total mass $M_0$, which indicates the critical mass for this particular initial condition lies between $38.8$ and $39$. 
 \rev{
	\paragraph{Dependence on physical and biological parameters.}
 Here we investigate the dependence of critical mass on other physical coefficients in the KS systems. Here we take the base configuration as specified in Section \ref{sec:agrregation} with physical coefficients set \eqref{eg_para}. We then change one of the coefficients and apply the same methodology as the aforementioned example, i.e. calculating $\frac{|c|_{\infty, H=24}}{|c|_{\infty, H=12}}$ for various values of the coefficient and find the first interval that $\frac{|c|_{\infty, H=24}}{|c|_{\infty, H=12}}$ is significantly away from $1$. In Tab. \ref{tab:coeff}, we summarize the results by examples that change one of the coefficients in \eqref{eg_para}. From Tab. \ref{tab:coeff}, we call tell that the following factors have a positive correlation with the critical mass (suppressing blow-up): mobility of the bacteria $\mu$ and chemical decay constant $k$. In contrast, the following factors have a negative correlation (promoting blow-up): chemo-sensitivity of the bacteria $\chi$ and time scale of chemotaxis $\epsilon$. 
 
 \begin{table}
     \centering
     \begin{tabular}{cccc}
     \toprule
         Original setup & New setup & new interval &  changes  \\
         \midrule
        $\mu=1$  & $\mu=0.8$ & $(38,38.2)$ & $\downarrow$  \\
         $\chi=1$ & $\chi=0.8$ & $(59.6,59.8)$ & $\uparrow$  \\
         $\epsilon=10^{-4}$ & $\epsilon=8\times 10^{-5}$ & $(47.8,48)$ &  $\uparrow$ \\ 
       $k=0.1$   & $k=0.08$ & $(46.8,47)$  & $\downarrow$  \\
       \bottomrule
     \end{tabular}
     \caption{Dependence of critical mass on KS physical and biological parameters. Original interval: $(47.6,47.8)$.}
     \label{tab:coeff}
 \end{table}
	}
	\subsection{Aggregation behaviors from non-radial initial data}
	In this subsection, we investigate aggregation behaviors in more general distributions. To this end, we \rev{present two experiments: 1) the initial distribution is not centered at the origin or even not a standard Fourier collocation point; 2) a more practical scenario where the initial distribution $\rho$ models several separated clusters of organisms. }

\rev{
 \paragraph{Shifted initial distribution.} In this example, we re-do the computation to find the critical mass in parameter setup \eqref{eg_para} as the same approach as Section \ref{sec:blowup_num}, while \emph{shifting the initial distribution concentrated at $[1/6,1/6,1/6]^T$}. The concentration point is selected away from the standard Fourier collocation point in our maximal resolution, namely $[-4,4]^3$ domain with each direction discretized by $H=24$ Fourier mode. Obviously, in this setup, the critical mass shall remain the same as the non-shifted one, namely $(47.6,47.8)$.

  In Fig. \ref{fig:eg_shift}(a), we compute the ration $\frac{|c|_{\infty, H=24}}{|c|_{\infty, H=12}}$ as in Fig.  \ref{fig:eg_mass}(b). The projected critical mass falls between $47.6$ and $47.8$ identical as the experiment result in Section \ref{sec:blowup_num}. In Fig. \ref{fig:eg_shift}(b) we show the location of particles representing $\rho$ at $T=1$. More specifically, we compute a histogram of particles projected to $x$-$y$ plane and zoom in to $[-1/3,1/3]^2$. Notice that $[-1/3,1/3]$ covers an interval of three Fourier collocation points under the resolution of computation. 
 
  Results in Fig. \ref{fig:eg_shift} confirm the capability of our method in predicting critical solutions that may focus on any location.
  In essence, our algorithm approximates $c$ with Fourier series  \eqref{trig_ser} and update of $c$ in Alg.\ref{alg:SIPM-C}, always shifting the $c$ series relative to the position of $X$ during collocation, see \eqref{eqn:sp_shft}. Hence our method avoids the potential numerical inaccuracy in \cite{fatkullin2012study} as pointed out in \cite{tomasevic:tel-01932777}.

 \begin{figure}[htbp]
		\centering
		
 		\subfigure[$\frac{|c|_{\infty, H=24}}{|c|_{\infty, H=12}}$ at $T=1$ vs. $M_0$. ]{\includegraphics[width=0.45\linewidth]{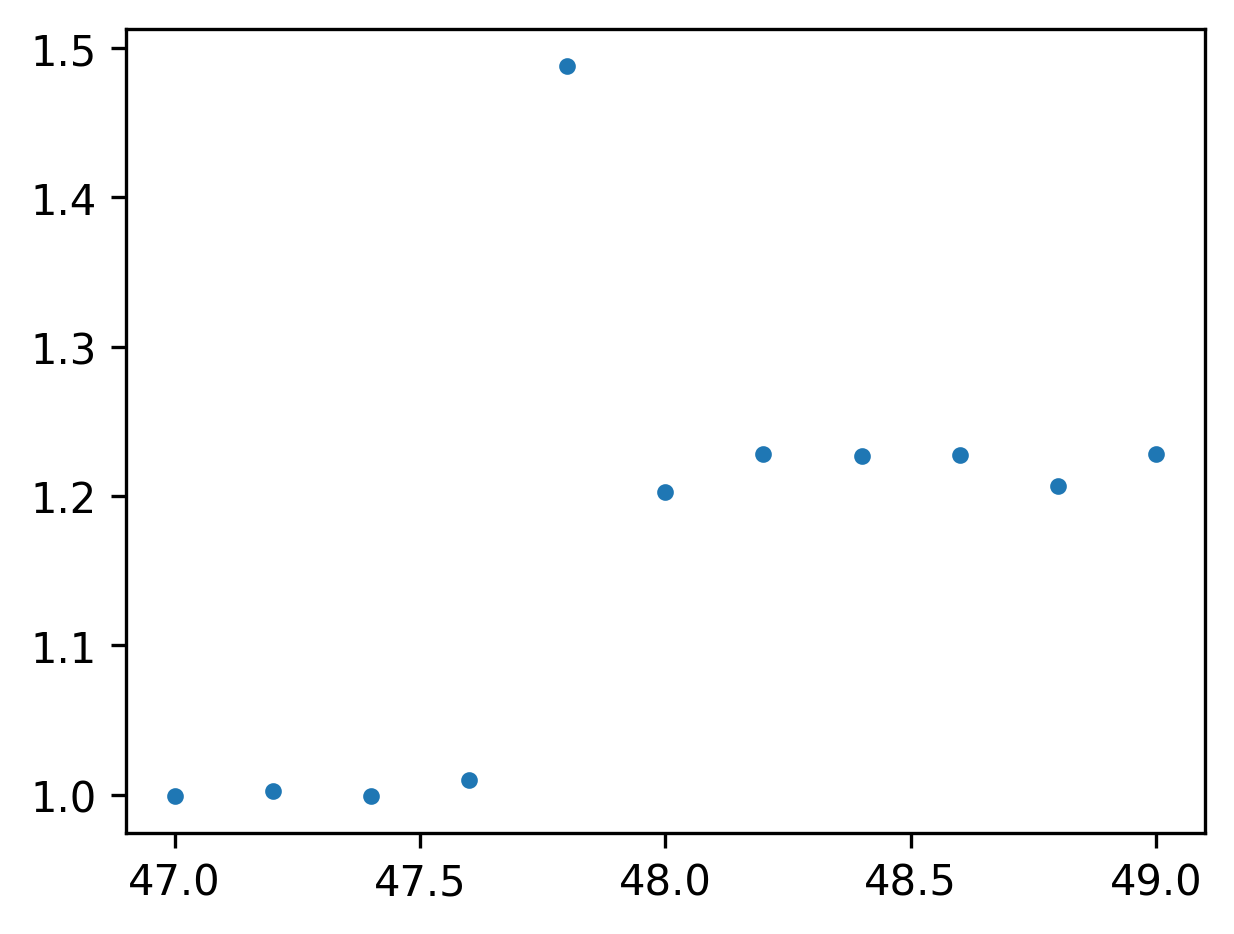}}
   \subfigure[Histogram of $X$ representing $\rho$ at $T=1$ with $M_0=47.8$, projected to $x-y$ plane.]{\includegraphics[width=0.45\linewidth]{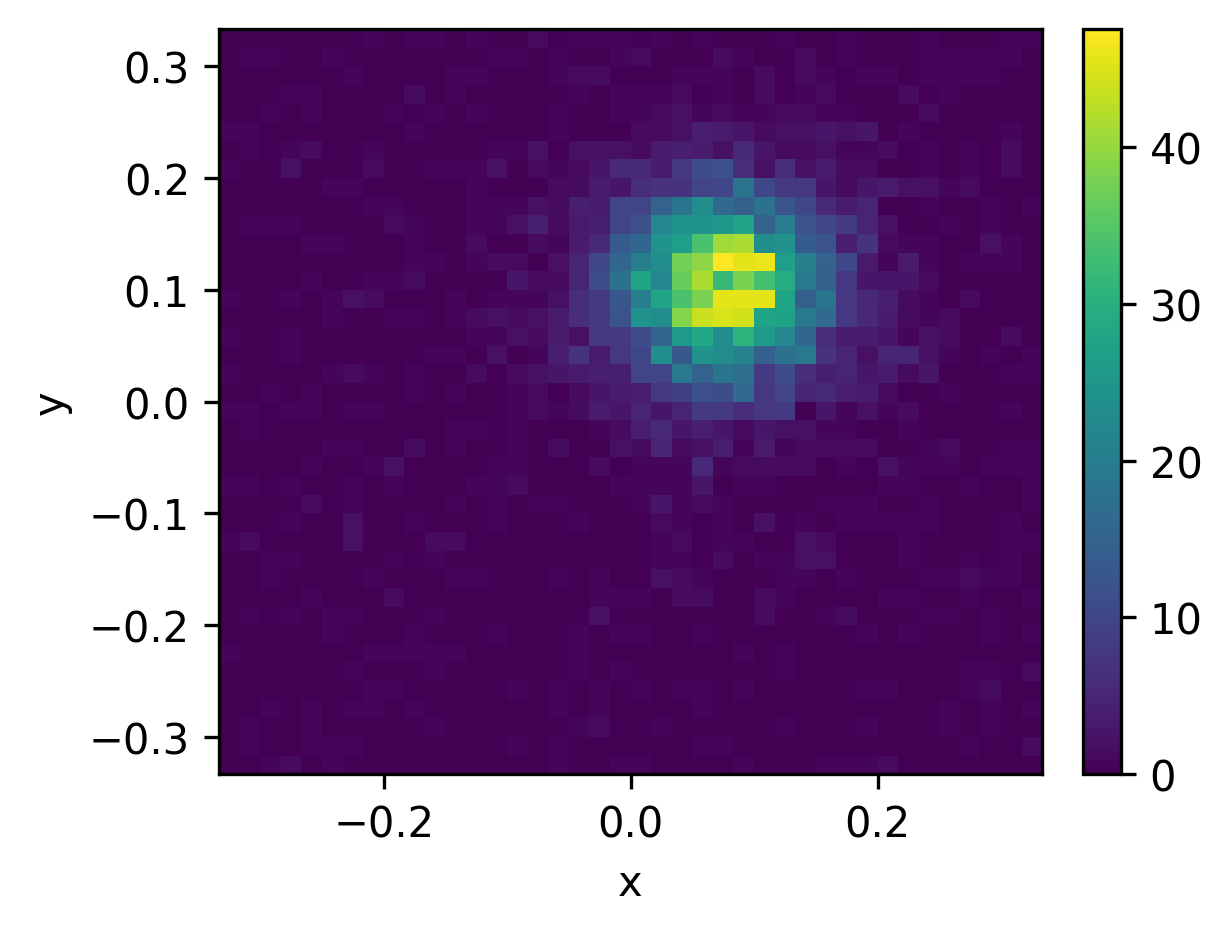}}
		
		\caption{Results in the example of shifted initial distribution: SIPF predicts blow up with concentration not on the regular Fourier grids.}
 		\label{fig:eg_shift}
 	\end{figure}
 }
 
\begin{figure}[htbp]
		\centering
		\subfigure[Initial distribution.]{\includegraphics[width=0.45\linewidth]{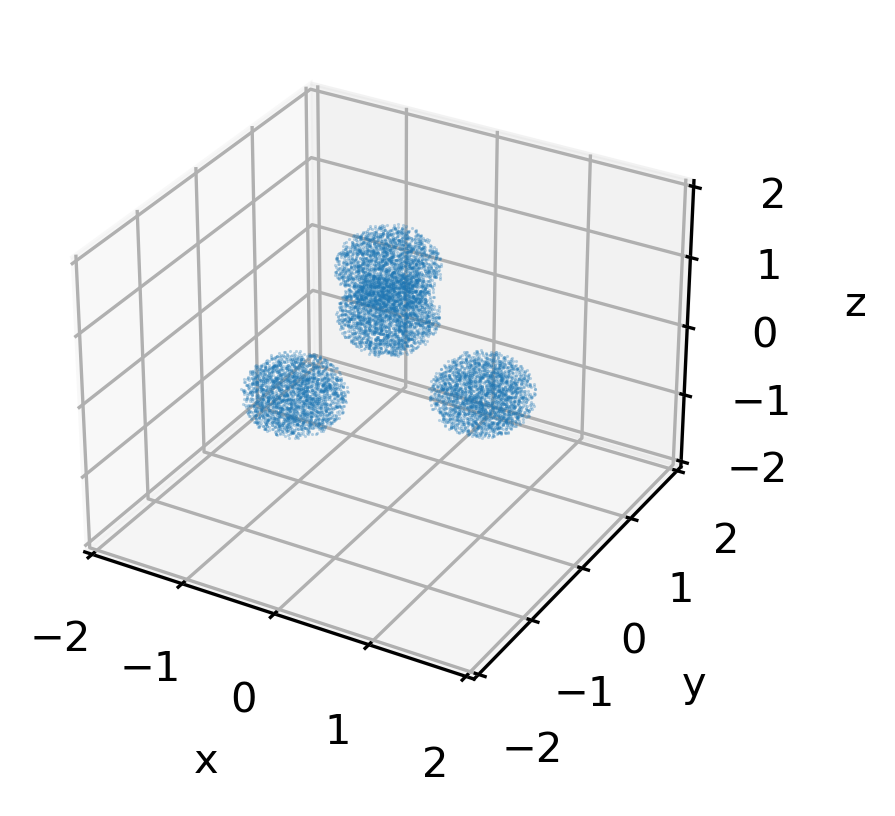}}
		\subfigure[$\frac{|c|_{\infty, H=24}}{|c|_{\infty, H=12}}$ vs. computation time.]{\includegraphics[width=0.45\linewidth]{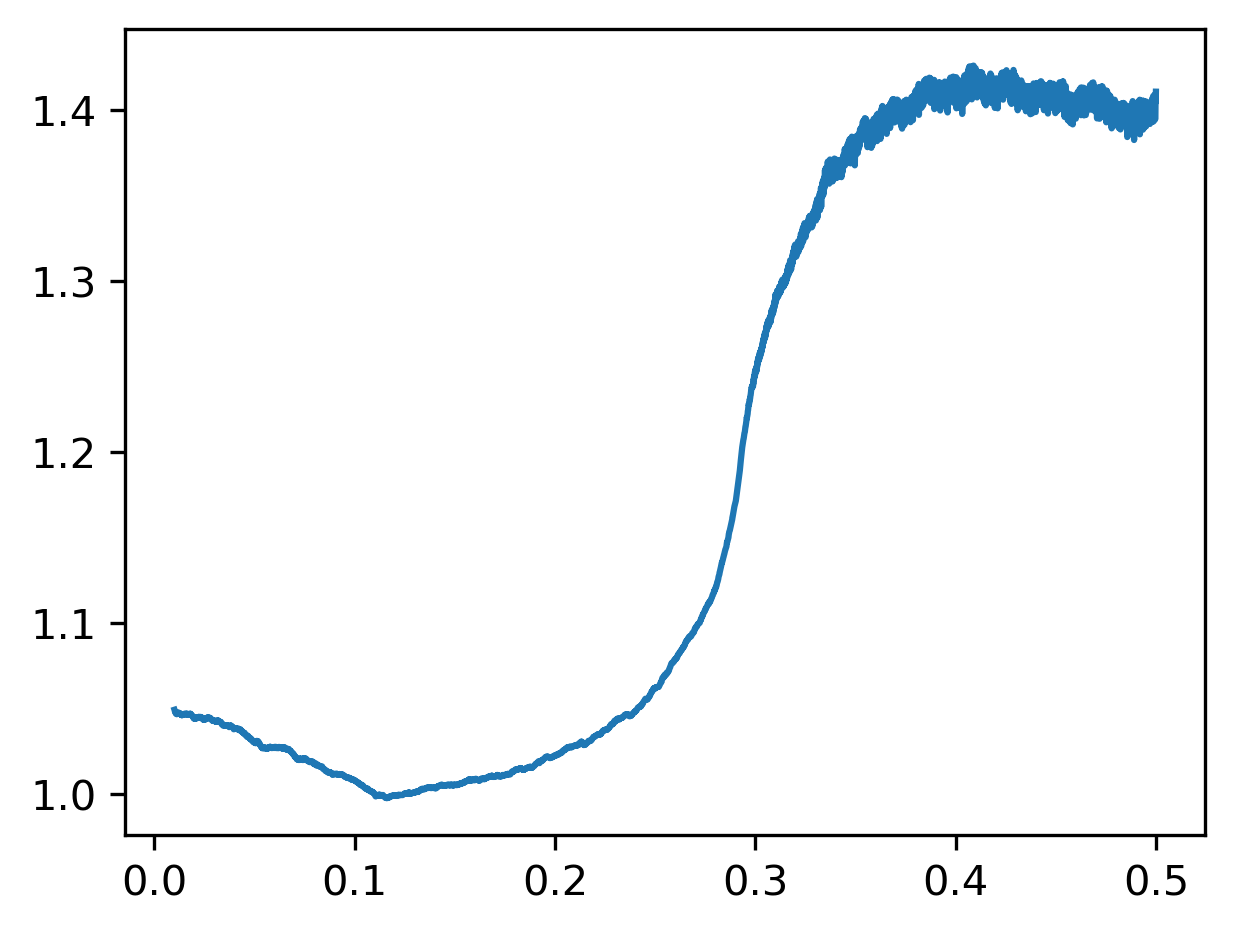}}
		\caption{Identifying the formation of a finite time singularity at $t\approx 0.3$ in non-radial solutions.}
		\label{fig:eg4_f1}
	\end{figure}

 \rev{\paragraph{Multi-clustered initial distribution.} This is a more practical scenario where the initial distribution $\rho$ models several separated clusters of organisms.} The mass in each cluster is below the critical mass, while the total mass is super-critical. To be more specific, we assume the initial distribution is a uniform distribution on four balls with a radius of $0.5$. These balls are centered at four vertices of a regular tetrahedron, namely,
	\begin{align}
		\begin{pmatrix}
			1\\0\\0
		\end{pmatrix}, \quad
		\begin{pmatrix}
			-\frac{1}{2}\\\frac{\sqrt{3}}{2}\\0
		\end{pmatrix}, \quad
		\begin{pmatrix}
			-\frac{1}{2}\\-\frac{\sqrt{3}}{2}\\0
		\end{pmatrix}, \quad
		\begin{pmatrix}
			0\\0\\\sqrt{2}
		\end{pmatrix}.
	\end{align}
	See also Fig. \ref{fig:eg4_f1}(a) for the scatter plot of particles representing the initial distribution.
	In this case, we assume the total mass is $M_0=80$ and each cluster has a mass of $20$, which is below the critical mass for a ball with a radius of $r=0.5$. Next, we apply the algorithm to compute the KS system up to $T=0.5$, with two different spatial discretizations ($H=24$ and $H=12$), while keeping the rest of the configurations. In Fig. \ref{fig:eg4_f1}(b), we calculate the ratio between the maxima of $c$ versus time for the two different spatial discretizations. We can see the singularities formed in the system at around $T=0.3$. 
	
\begin{figure}[htbp]
		\centering
		\subfigure[$T=0.1$]{\includegraphics[width=0.45\linewidth]{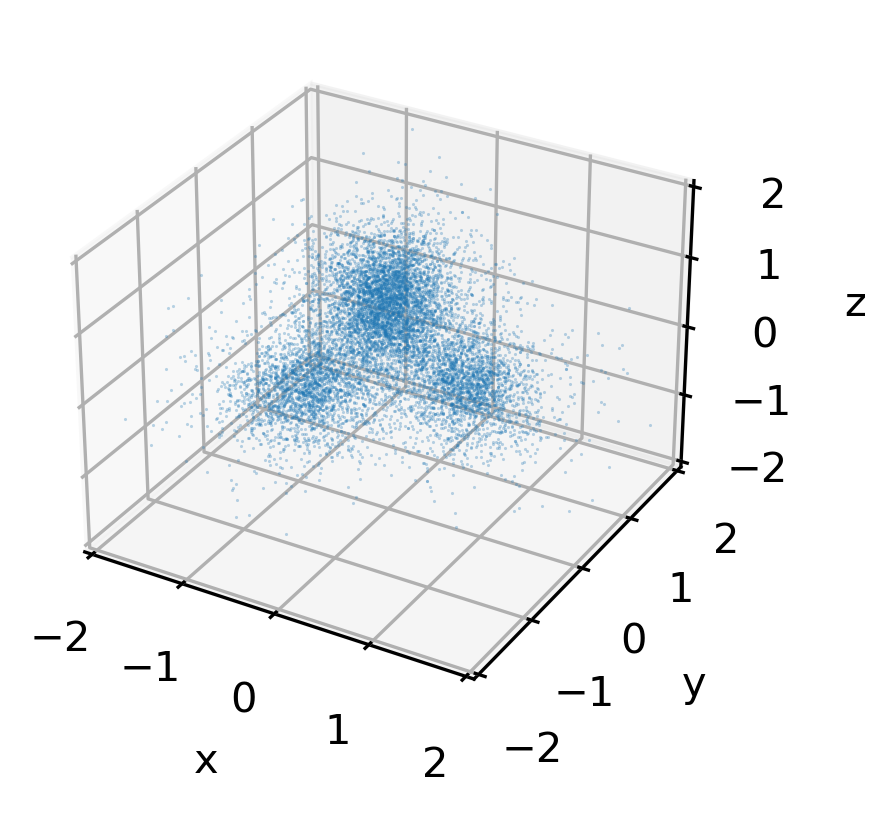}}
		\subfigure[$T=0.2$]{\includegraphics[width=0.45\linewidth]{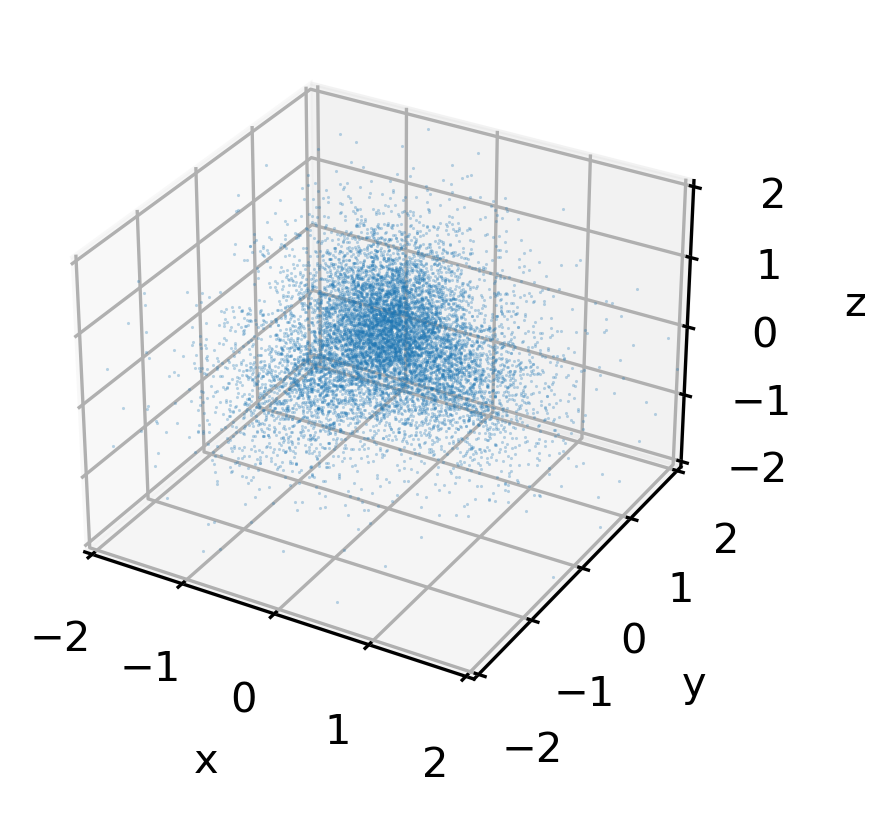}}\\
		\subfigure[$T=0.3$]{\includegraphics[width=0.45\linewidth]{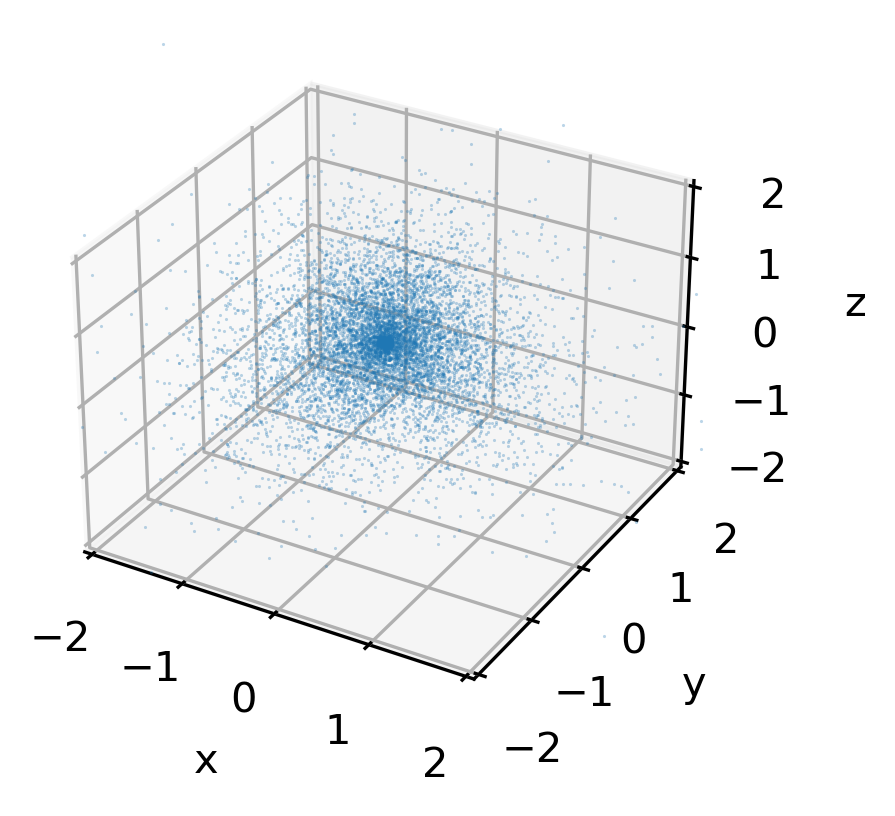}}
		\subfigure[$T=0.4$]{\includegraphics[width=0.45\linewidth]{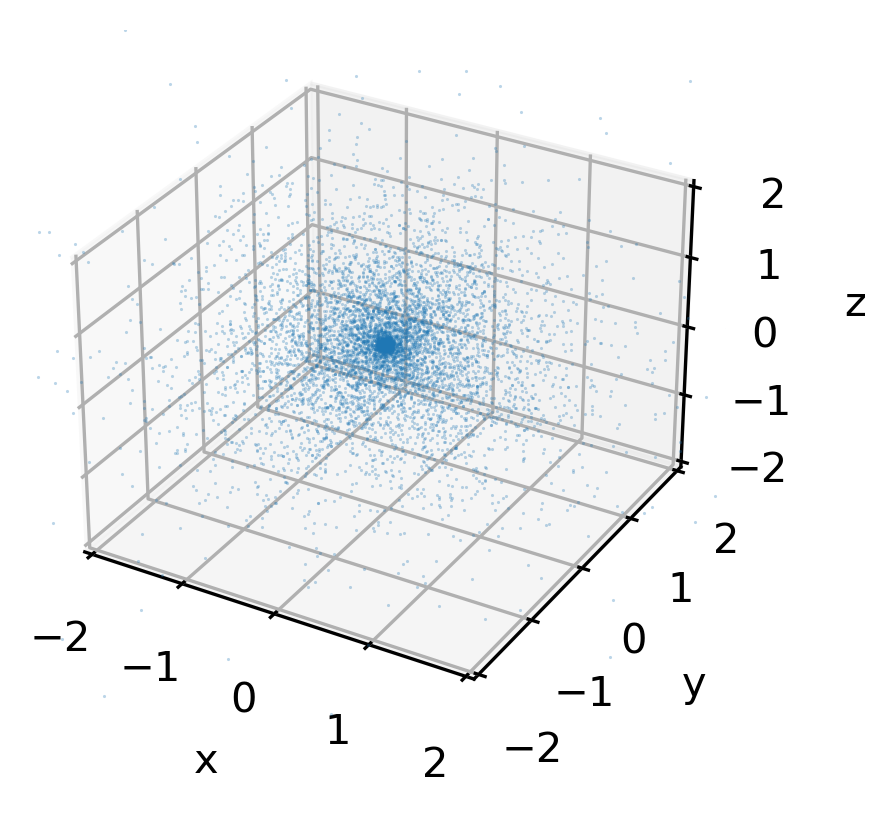}}
		\caption{Particle scatter plot at $T=0.1:0.1:0.4$: three cluster merging and a singularity formation.}
		\label{fig:eg4_f2}
	\end{figure}

	In Fig. \ref{fig:eg4_f2}, we present the scatter plot of particles between the time $T=0.1$ and $T=0.4$.
	By comparing Fig. \ref{fig:eg4_f1}(a) with Fig. \ref{fig:eg4_f2}(a), we can see diffusive behavior. This behavior is a result of the mass in each cluster being below the critical mass. The diffusive behavior persists until approximately $T=0.2$, as depicted in Fig. \ref{fig:eg4_f1}(b), where the active particles form a single larger cluster. The mass of this new cluster, centered at the origin, is $M_0=80$. In Fig. \ref{fig:eg4_f2}(c), we can observe the aggregation process 
	starting to form a singularity. This can also be seen from the sharp increase in the ratio of the maximum of $c$ in Fig. \ref{fig:eg4_f1}(b). Finally, in Fig. \ref{fig:eg4_f2}(d), we can directly identify the possible blow-up at the origin through the scatter plot.
	
	\subsection{Critical mass and blowup in  parabolic-parabolic KS}
	As the last example, we examine the singular solutions in the fully parabolic KS systems. For the purpose of exposition, we set $\epsilon=0.1$ in \eqref{ppKS}, while keeping the remaining physical parameters constant. The initial condition is assumed to be a uniform distribution on a ball with a radius of $0.8$ and $c(x,0)=0$. 
	
	From Fig. \ref{fig:eg_mass1}, we can determine that the critical mass is approximately $M_0=39$. We apply the same computational configuration as in Fig. \ref{fig:eg_mass1}, except we enlarge the domain to $L=12$ to accommodate possible diffusive behaviors. We test our algorithm for two cases, $M_0=40$ and $M_0=160$, respectively.
	
\begin{figure}[htbp]
		\centering
		\subfigure[Scatter plot of particles at $T=1$ with $M_0=40$.]{\includegraphics[width=0.32\linewidth]{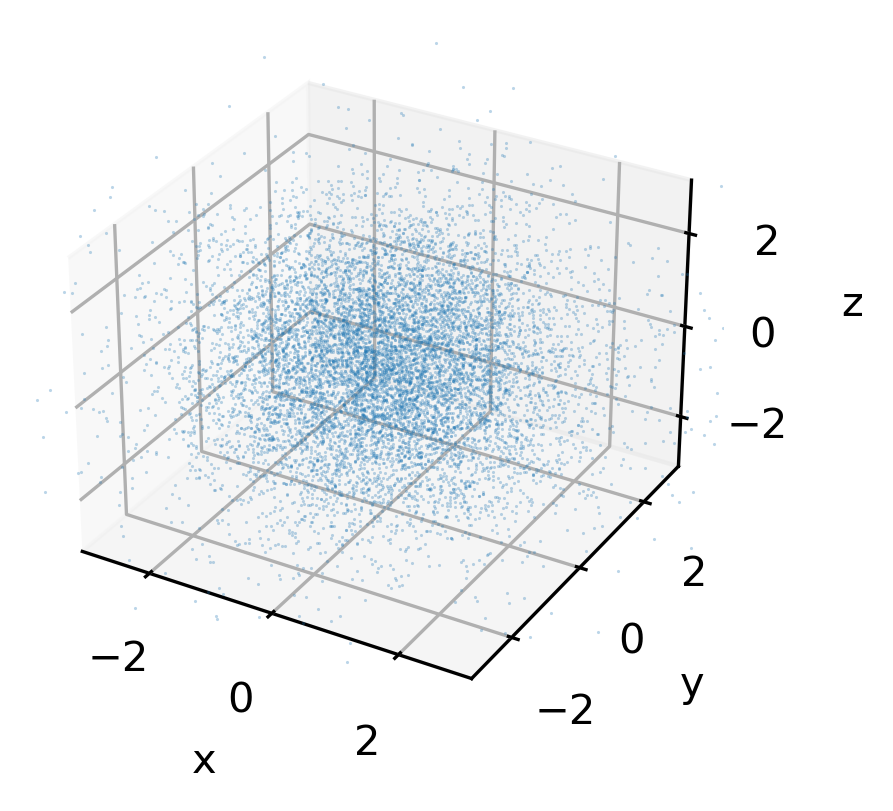}}
		\subfigure[Scatter plot of particles at $T=1$ with $M_0=160$.]{\includegraphics[width=0.32\linewidth]{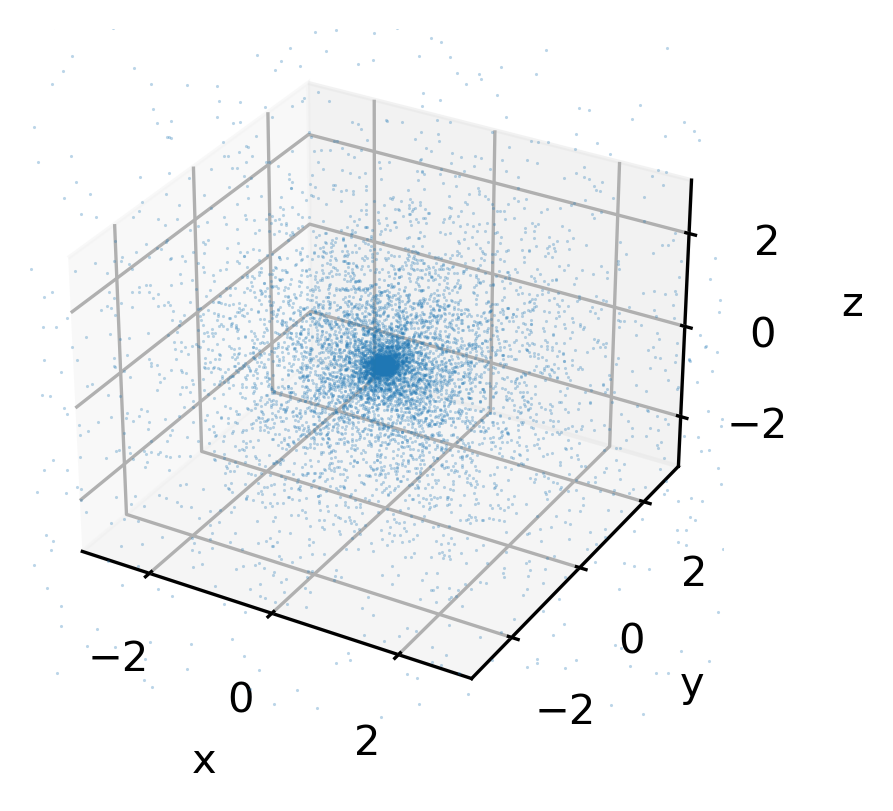}}
		\subfigure[$\frac{|c|_{\infty, H=24}}{|c|_{\infty, H=12}}$ vs computation time $T$ with different total mass $M_0$.]{\includegraphics[width=0.32\linewidth]{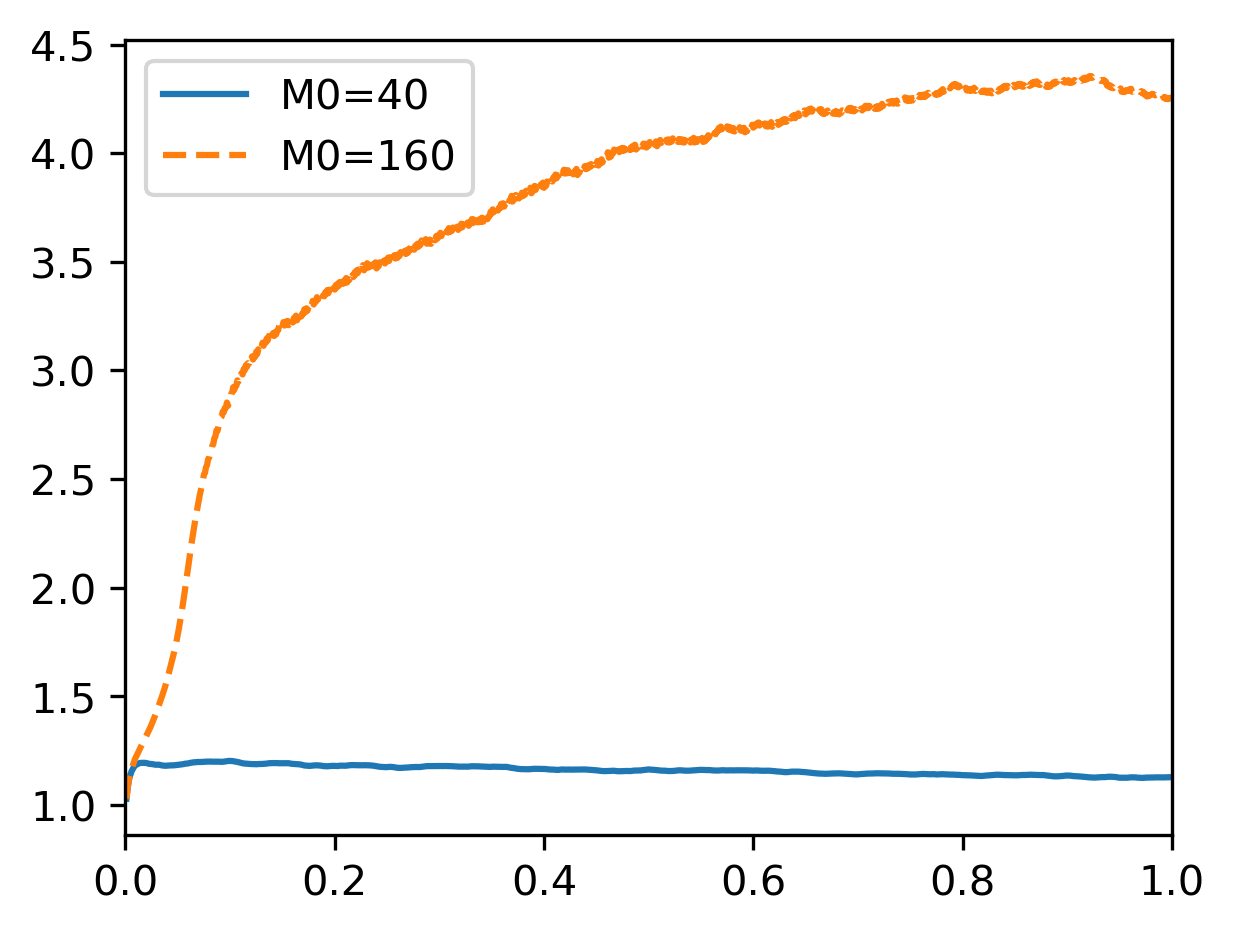}}
		
		\caption{Effects of initial mass $M_0$ on focusing behavior (finite time blowup).}
		\label{fig:eg_para}
	\end{figure}

	The behaviors of the system are reported in Fig. \ref{fig:eg_para}. In Fig. \ref{fig:eg_para}(a) and (b), we present the scatter plot of the particles representing the density $\rho$ with $M_0=40$ and $M_0=160$,  respectively. 
 
 We find that despite the initial mass $M_0=40$ being larger than the critical mass in the case of $\epsilon=10^{-4}$, the system does not blow up. We report that the variance of the particles grows linearly in computational time $T$, with diffusion coefficients fitted to be $1.696$. In the absence of the chemical attractant, namely $\chi=0$, the diffusion coefficient is expected to be $4\mu=4$. \rev{This implies that the parabolic-parabolic KS systems with mass below critical mass are effectively diffusive, with diffusion suppressed by chemical attraction.}
 
 However, for $M_0=160$, the system exhibits a possible singularity at the origin. In Fig. \ref{fig:eg_para}(c), we present the ratio of $|c|_\infty$ under $H=24$ and $H=12$ for both initial masses. Similar to the observation in Fig. \ref{fig:eg_para}(a) and (b), the blow-up behavior crucially depends on a critical level of the initial mass.
	\section{Concluding Remarks}
	In this paper, We developed a stochastic interacting particle and field algorithm, observed its convergence, and demonstrated its efficacy in computing blowup dynamics of fully parabolic KS systems in 3D from general non-radial initial data. The algorithm is recursive and does not have any history dependence, and the field variable is computed using \rev{Fourier series}. Since the field variable (concentration) is smoother than the density, the \rev{series} approach works well with only a few Fourier modes. The aggregation or focusing behavior in the density variable is resolved by using 10k particles. The algorithm successfully detects blowup through the field variable based on the critical amount of initial mass. The algorithm is self-adaptive and does not rely on any assumption about the blowup behavior, which is unknown except in the case of the parabolic-elliptic KS system. A potential weakness of the algorithm is the high cost of \rev{series expansion in 3D} when a large number of Fourier modes is required for high-resolution computation near the blowup time. We will study this issue in future work.

	\section*{Acknowledgements}
	\noindent
	ZW was partially supported by NTU SUG-023162-00001, MOE AcRF Tier 1 Grant RG17/24, and JX by NSF grant DMS-2309520. ZZ was supported by the National Natural Science Foundation of China  (Projects 12171406 and 92470103), the Hong Kong RGC grant (Projects 17307921 and 17304324), the Seed Funding Programme for Basic Research (HKU), the Outstanding Young Researcher Award of HKU (2020-21), and Seed Funding for Strategic Interdisciplinary Research Scheme 2021/22 (HKU). 

 \section*{Declarations}
The authors have no competing interests to declare that are relevant to the content of this article.

\bibliographystyle{plain}
	\bibliography{ZWpaper_DeepLearningPDEDiscCoef} 

\end{document}